\documentclass[11pt, a4paper]{amsart}
\usepackage{amsfonts}

\newtheorem{theorem}{\bf Theorem}[section]
\newtheorem{remark}{\bf Remark}[section]
\newtheorem{proposition}{Proposition}[section]

\newtheorem{corollary}{Corollary}[section]

\newtheorem{example}{Example}[section]
\newtheorem{definition}{Definition}[section]

\usepackage[utf8]{inputenc}
\usepackage{diffcoeff}

\usepackage{enumerate,mathrsfs}
\usepackage{amssymb,amsmath,graphicx,amsthm,mathtools,hyperref}

\usepackage{booktabs}
\usepackage{epstopdf}

\usepackage{subfig}

\theoremstyle{plain}
\usepackage{siunitx}

\usepackage{algorithm}
\usepackage{algorithmic,colortbl}

\makeindex
\usepackage[intoc]{nomencl} 

\makenomenclature

\usepackage{hyperref}
\hypersetup{nesting=true,debug=true,naturalnames=true}
\usepackage{graphicx,amssymb,upref}

\usepackage{enumerate,float}
\usepackage{fullpage}
\makeatletter
\@namedef{subjclassname@2020}{\textup{2020} Mathematics Subject Classification}
\makeatother

\usepackage{lipsum}
\begin{document}

\title{Fractional counting process at  L\'{e}vy times and its applications   }
	\author[Shilpa Garg]{S\lowercase{hilpa} G\lowercase{arg}$^{1*}$,~A\lowercase{shok} K\lowercase{umar} P\lowercase{athak}$^{2*}$,~A\lowercase{ditya} M\lowercase{aheshwari}$^{\#}$
\address{	$^*$D\lowercase{epartment of} M\lowercase{athematics and} S\lowercase{tatistics}, C\lowercase{entral} U\lowercase{niversity of} P\lowercase{unjab},
	B\lowercase{athinda}, P\lowercase{unjab}-151401, I\lowercase{ndia}.}
  \address{  $^{\#}$O\lowercase{perations} M\lowercase{anagement and} Q\lowercase{uantitative} T\lowercase{echniques} A\lowercase{rea}, I\lowercase{ndian} I\lowercase{stitute of} M\lowercase{anagement} I\lowercase{ndore}, I\lowercase{ndore}, M\lowercase{adhya} P\lowercase{radesh} - 453556, INDIA}\\
\email{\lowercase{shilpa17garg@gmail.com$^1$, ashokiitb09@gmail.com$^2$, adityam@iimidr.ac.in$^{\#}$}}}

	%E\lowercase{mail: ashokiitb09@gmail.com}
%	\\$^{2}$D\lowercase{epartment of} M\lowercase{athematics}, I\lowercase{ndian} I\lowercase{nstitute of} T\lowercase{echnology}
%	B\lowercase{ombay},\\ P\lowercase{owai}, M\lowercase{umbai}-400076, I\lowercase{ndia}.}
%E\lowercase{mail: pv@math.iitb.ac.in}}
%\thanks{Corresponding Author E-mail Address: ashokiitb09@gmail.com}%\\

%\thanks{The research of  R. Soni was supported by CSIR, Government of India.}

%\author[Ritik Soni]{R. Soni$^{1}$}
%\address{\emph{Department of Mathematics and Statistics,
%		Central University of Punjab, Bathinda, Punjab -151401, India.}}
%\email{ritiksoni2012@gmail.com}
%\author[]{Ashok Kumar pathak}
%\address{\emph{Department of Mathematics and Statistics,
%		Central University of Punjab, Bathinda, Punjab -151401, India.}}
%\email{ashokiitb09@gmail.com}

%\begin{center}
%{Neha Gupta}$^{\textrm{a}}$, {Aditya Maheshwari}$^{\textrm{b}}$
%\footnotesize{
	%		$$\begin{tabular}{l}
		%		$^{\textrm{a}}$ \emph{Department of Theoretical Statistics and Mathematics Unit,
			%Indian Statistical Institute Delhi, New Delhi - 110016, India.}
		%\end{tabular}$$
		%$$
		%\hspace{-1.2cm}\begin{tabular}{l}
			%		$^{\textrm{b}}$ \emph{Operations Management and Quantitative Techniques Area, Indian Institute of Management Indore, Indore 453556, India.}\\

			%\end{tabular}$$}
			%\end{center}
			%\vtwo
			%\begin{center}
			%\noindent{\bf Abstract}
			%\end{center}

			\begin{abstract} 
            Traditionally, fractional counting processes, such as the fractional Poisson process, \textit{etc} have been defined using fractional differential and integral operators. Recently, Laskin (2024) introduced a generalized fractional counting process (FCP) by changing the probability mass function (pmf) of the time fractional Poisson process using the generalized three-parameter Mittag-Leffler function. Here, we study some additional properties for the FCP and introduce a time-changed fractional counting process (TCFCP), defined by time-changing the FCP with an independent L\'{e}vy subordinator.  We derive distributional properties such as the  Laplace transform,  probability generating function, the moments generating function, mean, and variance for the TCFCP. Some results related to waiting time distribution and the first passage time distribution are also discussed. We define the multiplicative and additive compound variants for the FCP and the TCFCP and examine their distributional characteristics with some typical examples. We explore some interesting connections of the TCFCP with Bell polynomials by introducing subordinated generalized fractional Bell polynomials.  
            It is shown that the moments of the TCFCP can be represented in terms of the subordinated generalized fractional Bell polynomials. Finally, we present the application of the FCP in a shock deterioration model.

				\end{abstract} 
			\maketitle

			\section{Introduction}
		 \noindent The counting processes are extensively utilized stochastic process for studying the occurrence of random events across time. The well-known examples of counting process are the Poisson process, fractional Poisson process (FPP), and stretched Poisson process, which have numerous applications in traffic flow, queuing theory, telecommunications,  insurance, reliability theory, and physics (see \cite{guler2022forecasting}, \cite{jung1969risk},  \cite{laskin2009some}, \cite{stanislavsky2008two}, \cite{timmermann1999multiscale} and \cite{wang2022explicit}) due to its mathematical ease. Recently, Laskin \cite{laskin2024new} introduced a wider class of the fractional counting process (FCP) by directly constructing its probability distribution function. It is a generalization of prominent stochastic processes such as the Poisson process, FPP, and stretched Poisson process that have been studied in the recent past. Further, he demonstrated various applications of the FCP including fractional compound Poisson process, generalized fractional Bell polynomial, and introduced stretched quantum coherent states.
         
		 Recently, there has been significant interest in subordination based generalizations of the counting processes driven by both a theoretical and practical view.  The subordination techniques help in overcoming various limitations of the counting processes, for example, in the case of the Poisson process, it allows for non-exponential inter-arrival times, making the system more suitable for real-world phenomena.  Meerschaert et al. \cite{meerschaert2011fractional}  demonstrated that the FPP is a generalized form of the Poisson process, where the time variable is substituted by an independent inverse stable subordinator. Kumar et al. \cite{kumar2011time} explored the time-changed variant of the Poisson processes and studied its connection to fractional differential equations. For recent developments in this regard, one may consult to \cite{gajda2019stable, kumar2019linnik, maheshwari2023tempered, maheshwari2019fractional, orsingher2012space, soni2024generalizedspl}.		

This paper first considers a broader class of FCP studied in \cite{laskin2024new} and explores its additional properties and some results related to convergence. It is shown that the FCP does not possess independent increment property. We have discussed the order statistic associated with the FCP and derived its conditional distribution.  Next, we define an FCP at  L\'{e}vy times,, namely, time-changed fractional counting processes (TCFCP), which is obtained by time-changing the FCP with an independent L\'{e}vy subordinator. We derive its distributional properties and calculate its mean and variance formula. The waiting times and first passage time of the TCFCP have been evaluated. Additionally, we introduce a multiplicative compound fractional counting process (MCFCP), which can be viewed as a multiplicative compound variant of the FCP. We calculate its cumulative distribution function (cdf), and thereafter, we discuss the MCFCP at L\'{e}vy times. We also present some typical examples. 

We extend the fractional generalized compound process (FGCP) introduced by \cite{laskin2024new} by discussing its distributional properties.  Further, we introduce the FGCP at L\'{e}vy times and explore its important properties with examples. Next, we introduce the subordinated version of generalized fractional Bell polynomials studied in  \cite{laskin2024new}, which is calculated using the probability distribution function of the TCFCP. For the case of $\alpha$-stable subordinator, the subordinated generalized fractional Bell polynomials (SGFBP) correspond to the generalized fractional Bell polynomials.  A shock deterioration model governed by the FCP in which a random variable with a gamma distribution is used to model the magnitude of each shock in order to explain the shock deterioration of structural resistance is proposed and studied. It is known that (see \cite{biard2014fractional, maheshwari2016long,soni2024generalizedspl}) the fractional versions of the counting process introduce memory effects in arrival occurrences, and it becomes more useful in modelling (see \cite{bapat2024modelling, kumar2020fractional,karagiannis2004long}) real-life scenarios. Our proposed model leverages the same memory effect, and we believe it will be very practical in real-life applications.  

 In recent literature, counting processes in both one-dimensional and multidimensional forms, including time-changed variants, are commonly used in reliability theory, specifically in shock and deterioration models.  Wang \cite{wang2022explicit} introduced a compound Poisson process-based shock deterioration model for the reliability assessment of aging structures.  Di Crescenzo and Meoli \cite{di2023competing} considered the shock model governed by the bivariate Poisson process. Further,  Soni et al. \cite{soni2024bivariate}  extended the work of Di Crescenzo and Meoli \cite{di2023competing} by introducing a bivariate tempered model and studied competing risks and shock model governed by it. For more insights, one may consult \cite{pellerey1994shock,crowder2001classical,cha2016new}. In order to appropriately describe a system where past shocks influence future occurrences, the FCP takes memory effects into account. Furthermore, the FCP permits greater flexibility in interarrival timings, which is crucial in portraying the variable nature and heavy-tailed character of real-world shock.
 
    The structure of the paper is as follows: In Section \ref{sec:prelims}, we present some preliminary notations and definitions. In Section \ref{sec:some-res}, we discuss additional properties and convergence results of the FCP. In Section \ref{sec:tcfcp}, we introduce the time-changed version of the FCP by an L\'{e}vy subordinator, namely, TCFCP, and discuss its distributional properties along with mean and variance. The waiting times and first passage time of the TCFCP are also evaluated. We discuss the MCFCP and 
 the MCFCP at L\'{e}vy times, along with their distributional properties and some special cases. Next, the FGCP at L\'{e}vy times are presented with examples.  In Section \ref{sec:appl}, the SGFBP are introduced and their generating function is evaluated. Finally, we also present an application of the FCP in a shock deterioration model.

\section{Preliminaries}\label{sec:prelims}
In this section, we enlist some basic definitions and notations that will be used in subsequent sections. Let $\mathbb{N}$ denotes the set of natural numbers and $\mathbb{N}_0 = \mathbb{N} \cup \{0\}$. Let $\mathbb{R}$ denote the set of real numbers.
\subsection{L\'{e}vy Subordinator} A L\'{e}vy subordinator denoted by $H (t),\;t \geq 0$ is a  non-decreasing L\'{e}vy process with Laplace transform (LT)  (see \cite[Section 1.3.2]{applebaum2009levy})
  \begin{equation*}
  \mathbb{E}\left[e^{-sH(t)}\right] = e^{-t \psi(s)},\;\; s \geq 0,
  \end{equation*}
  where $\psi(s)$ is Laplace exponent given by 
  \begin{equation*}
  \psi(s) = \eta s +\int_{0}^{\infty} (1-e^{-sx})  \nu (\mathrm{d}x), \;\; \eta \geq 0.
  \end{equation*}
  Here $\eta$ is the drift coefficient and $\nu$ is a non-negative L\'{e}vy measure on the positive half-line satisfying
  \begin{equation*}
  \int_{0}^{\infty} \min \{x,1\} \nu (\mathrm{d}x) < \infty, \;\;\; \text{ and } \;\;\; \nu ([0,\infty)) =\infty,
  \end{equation*}
  so that $H (t),\;t \geq 0$ has strictly increasing sample paths almost surely (a.s.). The probability density function of $H(t)$ is denoted by $h(x,t)$.
\subsection{Fractional counting process}
A wider class of fractional counting process (FCP) is introduced by \cite{laskin2024new} denoted by $\mathcal{N}_{\mu, \vartheta }^{\zeta ,\theta}(t),\; t\geq 0$ having the probability distribution function as
\begin{equation*}
\mathcal{P}_{\mu, \vartheta }^{\zeta ,\theta}(n,t)= (\zeta )_n \mathit{\Gamma } ( \vartheta )\frac{(\lambda_{\theta}t^{\theta})^n}{n!}\sum_{k=0}^{\infty}\frac{(\zeta +n)_{k}(-\lambda_{\theta}t^{\theta})^k}{k!\mathit{\Gamma } (\mu(n+k)+ \vartheta )},\;\;\;\;t\geq0,
\end{equation*}
where $n=0,1,2, ... , $ and the parameters $\mu, \vartheta ,\zeta ,\theta$ and $\lambda_{\theta}$ satisfy the following conditions $0<\mu\leq 1,\;\zeta >0,\; \vartheta \geq\mu\zeta ,\;0<\theta\leq 1,\;\lambda_{\theta}>0$.\\
The probability distribution function can also be represented as
\begin{align}\label{process}
\mathcal{P}_{\mu, \vartheta }^{\zeta ,\theta}(n,t)&= \mathit{\Gamma } ( \vartheta )\frac{(-z)^n}{n!}\frac{\mathrm{d}^n }{\mathrm{d} z^n}E_{\mu, \vartheta }^{\zeta }(z){\Big|}_{z=-\lambda_{\theta}t^{\theta}},\;\;\;\;t\geq 0,\;\;n=0,1,2,...,\\
\mathcal{P}_{\mu, \vartheta }^{\zeta ,\theta}(0,t)&= \mathit{\Gamma } ( \vartheta )E_{\mu, \vartheta }^{\zeta }(-\lambda_{\theta}t^{\theta}),\nonumber
\end{align}
where $E_{\mu, \vartheta }^{\zeta }(z)$ represents the generalized three-parameter Mittag-Leffler function given by (see \cite{prabhakar1971singular})
\begin{align*}
    E_{\mu, \vartheta }^{\zeta }(z) = \sum_{m=0}^{\infty}\frac{(\zeta )_m}{m!\mathit{\Gamma } (\mu m+ \vartheta )}z^m,
\end{align*}
where $(\zeta )_m$ and $\mathit{\Gamma }(\alpha)$ represent the Pochhammer symbol and Gamma function, respectively and are defined as
\begin{equation*}
    (\zeta )_m=\frac{\mathit{\Gamma } (\zeta  + m)}{\mathit{\Gamma } (\zeta )},\;\;\;\; \mathit{\Gamma }(\alpha)= \int_{0}^{\infty}e^{-t}t^{\alpha-1}\mathrm{d} t.
\end{equation*}
The derivative of $E_{\mu, \vartheta }^{\zeta }(z)$ of order $n$, for any $n\in\mathbb{N}$ is given by
\begin{align*}
    E_{\mu, \vartheta }^{\zeta (n)}(x) &= \frac{\mathrm{d} ^n}{\mathrm{d} z^n} E_{\mu, \vartheta }^{\zeta }(z)\big{|}_{z=x},
    \end{align*}
    \begin{align}\label{derivative}
    \frac{\mathrm{d} ^n}{\mathrm{d} z^n} E_{\mu, \vartheta }^{\zeta }(z) &= (\zeta )_n E_{\mu, \vartheta +n\mu}^{\zeta +n}(z).
    \end{align}
The probability generating function (pgf) of $\mathcal{N}_{\mu, \vartheta }^{\zeta ,\theta}(t),\; t\geq 0$ is given as 
\begin{equation}\label{pgf}
\mathcal{G}_{\mu, \vartheta }^{\zeta ,\theta}(s,t)=\mathit{\Gamma } ( \vartheta )E_{\mu, \vartheta }^{\zeta }(\lambda_{\theta}t^{\theta}(s-1)).
\end{equation} 
The moment generating function (mgf) of $\mathcal{N}_{\mu, \vartheta }^{\zeta ,\theta}(t),\; t\geq 0$ is given as \begin{equation}\label{mgf}
\mathcal{H}_{\mu, \vartheta }^{\zeta ,\theta}(s,t)=\mathit{\Gamma } ( \vartheta )E_{\mu, \vartheta }^{\zeta }(\lambda_{\theta}t^{\theta}(e^{-s}-1)).
\end{equation} 
 The mean and variance of the process defined in Eq.(\ref{process}) is given by
\begin{align}\label{mean}
\mathbb{E}\left[\mathcal{N}_{\mu, \vartheta }^{\zeta ,\theta}(t)\right]&= \frac{\zeta  \mathit{\Gamma } ( \vartheta )\lambda_{\theta}t^{\theta}}{\mathit{\Gamma }(\mu+ \vartheta )},\\
\label{var}
\text{Var}\left[\mathcal{N}_{\mu, \vartheta }^{\zeta ,\theta}(t)\right]&= \mathbb{E}\left[\mathcal{N}_{\mu, \vartheta }^{\zeta ,\theta}(t)\right] +\left [\mathbb{E}\left[\mathcal{N}_{\mu, \vartheta }^{\zeta ,\theta}(t)\right] \right ]^{2}\left\{ \left(1+\frac{1}{\zeta }\right)\frac{B(\mu+ \vartheta ,\mu+ \vartheta )}{B(2\mu+ \vartheta , \vartheta )}-1\right\},
\end{align}
where $B(\alpha,\beta)$ is the beta function defined as
\begin{equation*}
   B(\alpha,\beta)= \frac{\mathit{\Gamma } (\alpha)\mathit{\Gamma } (\beta)}{\mathit{\Gamma } (\alpha+\beta)}.
\end{equation*}
\section{Some results of the FCP}\label{sec:some-res}
In this section, we present some additional properties and convergence results for the FCP.
Here, we first determine the LT of the FCP as defined by Eq. (\ref{process}). The LT, denoted by $\mathcal{L}_{\mu, \vartheta }^{\zeta ,\theta}(s,t)$, can be evaluated using the pgf provided in Eq. (\ref{pgf})
	\begin{align}\label{LT}
		\mathcal{L}_{\mu, \vartheta }^{\zeta ,\theta}(s,t)&=\mathit{\Gamma } ( \vartheta )E_{\mu, \vartheta }^{\zeta }(\lambda_{\theta}t^{\theta}(e^{-s}-1))\\
		&= \mathit{\Gamma } ( \vartheta )\sum_{m=0}^{\infty}\frac{(\zeta )_m(\lambda_{\theta}t^{\theta}(e^{-s}-1))^m}{m!\mathit{\Gamma } (\mu m+ \vartheta )}.\nonumber
	\end{align} 

\subsection{Lack of independent increments property}
The FCP loses the independent increment property and, therefore, does not possess the lack of memory property. It is explained with the help of LT as given in Eq. (\ref{LT})
\begin{align*}
	\mathcal{L}_{\mu, \vartheta }^{\zeta ,\theta}(s,t_1+t_2)&=\mathit{\Gamma } ( \vartheta )E_{\mu, \vartheta }^{\zeta }(\lambda_{\theta}(t_1+t_2)^{\theta}(e^{-s}-1))\\
	&\neq \mathit{\Gamma } ( \vartheta )E_{\mu, \vartheta }^{\zeta }(\lambda_{\theta}t_1^{\theta}(e^{-s}-1))\mathit{\Gamma } ( \vartheta )E_{\mu, \vartheta }^{\zeta }(\lambda_{\theta}t_2^{\theta}(e^{-s}-1)).
\end{align*}
Therefore,  we get
\begin{align*}
	\mathcal{L}_{\mu, \vartheta }^{\zeta ,\theta}(s,t_1+t_2) \neq \mathcal{L}_{\mu, \vartheta }^{\zeta ,\theta}(s,t_1) \mathcal{L}_{\mu, \vartheta }^{\zeta ,\theta}(s,t_2).
\end{align*}
%Hence, the FCP has no independent increments.

\subsection{Order statistics}
Let $Y_1,Y_2,\cdots,Y_n$ be $n$ iid random variables with pdf $g_Y$. Consider  $Y_{(1)}\leq Y_{(2)}\leq\cdots \leq Y_{(k)}\cdots \leq Y_{(n)}$, then 
$(Y_{(1)},Y_{(2)},\cdots,Y_{(n)})$ is the order statistics. Let $ k\in \{1,2,\cdots,\mathcal{N}_{\mu, \vartheta }^{\zeta ,\theta}(t)\}$ and we denote   $Y_{(k)}^{\mathcal{N}_{\mu, \vartheta }^{\zeta ,\theta}(t)}$ as the $k$-th order statistics for $\mathcal{N}_{\mu, \vartheta }^{\zeta ,\theta}(t)$ samples.
\begin{theorem}
	Let $\mathcal{N}_{\mu, \vartheta }^{\zeta ,\theta}(t)$ be the FCP and $\{Y_i, i\in \mathbb{N}\}$ be iid random variables having probability distribution function $G_Y$, 
	\begin{equation*}
\mathbb{P}\left[Y_{(k)}^{\mathcal{N}_{\mu, \vartheta }^{\zeta ,\theta}(t)}<x \;\big{|} \; \mathcal{N}_{\mu, \vartheta }^{\zeta ,\theta}(t)\geq k\right]= \frac{\mathbb{P}\left[\tilde{\mathcal{N}}_{\mu, \vartheta }^{\zeta ,\theta}(t)\geq k\right]}{\mathbb{P}\left[\mathcal{N}_{\mu, \vartheta }^{\zeta ,\theta}(t)\geq k\right]}, \;\;\;\; k\in \mathbb{N},
	\end{equation*}
	where $\tilde{\mathcal{N}}_{\mu, \vartheta }^{\zeta ,\theta}(t)$ is the FCP with parameter $\lambda_{\theta}G_Y >0$.
\end{theorem}
\begin{proof}
	With the help of the law of conditional probability, we have that
	\begin{align}
		\mathbb{P}\left[Y_{(k)}^{\mathcal{N}_{\mu, \vartheta }^{\zeta ,\theta}(t)}<x\;\big{|} \;\mathcal{N}_{\mu, \vartheta }^{\zeta ,\theta}(t)\geq k\right]&= \frac{\sum_{n=k}^{\infty}\mathbb{P}\left[Y_{(k)}^{\mathcal{N}_{\mu, \vartheta }^{\zeta ,\theta}(t)}<x \;,\;\mathcal{N}_{\mu, \vartheta }^{\zeta ,\theta}(t)=n\right]}{\mathbb{P}\left[\mathcal{N}_{\mu, \vartheta }^{\zeta ,\theta}(t)\geq k\right]}\nonumber\\
		&= \frac{\sum_{n=k}^{\infty}\mathbb{P}\left[Y_{(k)}^{\mathcal{N}_{\mu, \vartheta }^{\zeta ,\theta}(t)}<x\;\big{|}\;\mathcal{N}_{\mu, \vartheta }^{\zeta ,\theta}(t)=n\right]\mathbb{P}\left[\mathcal{N}_{\mu, \vartheta }^{\zeta ,\theta}(t)=n\right]}{\mathbb{P}\left[\mathcal{N}_{\mu, \vartheta }^{\zeta ,\theta}(t)\geq k\right]}.\nonumber
	\end{align}
	Let us first simplify the numerator 
	\begin{align}
		\sum_{n=k}^{\infty}\mathbb{P}&\left[Y_{(k)}^{\mathcal{N}_{\mu, \vartheta }^{\zeta ,\theta}(t)}<x\;\big{|}\;\mathcal{N}_{\mu, \vartheta }^{\zeta ,\theta}(t)=n\right]\mathbb{P}\left[\mathcal{N}_{\mu, \vartheta }^{\zeta ,\theta}(t)=n\right]\nonumber\\ 
		=& \sum_{n=k}^{\infty}\sum_{j=k}^{n} \binom{n}{j}G_Y^j(x)(1-G_Y(x))^{n-j}(\zeta )_n \mathit{\Gamma } ( \vartheta )\frac{(\lambda_{\theta}t^{\theta})^n}{n!}\sum_{m=0}^{\infty}\frac{(\zeta +n)_{m}(-\lambda_{\theta}t^{\theta})^m}{m!\mathit{\Gamma } (\mu(n+m)+ \vartheta )}\nonumber\\
		=&\sum_{j=k}^{\infty}\sum_{m=0}^{\infty}\sum_{n=j}^{\infty} \frac{G_Y^j(x)(1-G_Y(x))^{n-j}}{(n-j)!j!}\frac{(\zeta )_n \mathit{\Gamma } ( \vartheta )(\lambda_{\theta}t^{\theta})^n(\zeta +n)_{m}(-\lambda_{\theta}t^{\theta})^m}{m!\mathit{\Gamma } (\mu(n+m)+ \vartheta )}\nonumber\\
		=&\sum_{j=k}^{\infty}\frac{G_Y^j(x)}{j!}\sum_{m=0}^{\infty}\frac{(-\lambda_{\theta}t^{\theta})^m \mathit{\Gamma } ( \vartheta )(\lambda_{\theta}t^{\theta})^j }{m!}\sum_{n=0}^{\infty}\frac{(1-G_Y(x))^{n}(\lambda_{\theta}t^{\theta})^n(\zeta )_{n+j}(\zeta+n+j )_{m}}{n!\mathit{\Gamma } (\mu(n+m+j)+ \vartheta )}\nonumber\\
		=&\sum_{j=k}^{\infty}\frac{G_Y^j(x)}{j!}\sum_{m=0}^{\infty}\frac{(-\lambda_{\theta}t^{\theta})^m \mathit{\Gamma } ( \vartheta )(\lambda_{\theta}t^{\theta})^j }{m!}\sum_{n=0}^{\infty}\frac{(1-G_Y(x))^{n}(\lambda_{\theta}t^{\theta})^n(\zeta )_{j}(\zeta+j )_{m}(\zeta+m+j )_{n}}{n!\mathit{\Gamma } (\mu(n+m+j)+ \vartheta )}.\label{thm-osnum}
	\end{align}
	The last step follows from $(\zeta )_{n+j}(\zeta+n+j )_{m}=(\zeta )_{j}(\zeta+j )_{m}(\zeta+m+j )_{n}$. Using Eq. (\ref{derivative}) in the above equation, we obtain
	\begin{align*}
		(\zeta+j )_{m} \sum_{n=0}^{\infty}\frac{(1-G_Y(x))^{n}(\lambda_{\theta}t^{\theta})^n(\zeta+m+j )_{n}}{n!\mathit{\Gamma } (\mu(n+m+j)+ \vartheta )}=&\;  (\zeta+j )_{m}E_{\mu, \vartheta +(m+j)\mu}^{\zeta +j}(\lambda_{\theta}(1-G_Y(x))t^{\theta})\\
		=&\; E_{\mu, \vartheta+\mu j }^{\zeta+j (m)}(\lambda_{\theta}(1-G_Y(x))t^{\theta}).
	\end{align*}
	Substituting the above in  Eq. (\ref{thm-osnum}) and with the help of Taylor's series, we get
	\begin{align*}
		\sum_{n=k}^{\infty}\mathbb{P}\left[Y_{(k)}^{\mathcal{N}_{\mu, \vartheta }^{\zeta ,\theta}(t)}<x\;\big{|}\;\mathcal{N}_{\mu, \vartheta }^{\zeta ,\theta}(t)=n\right]&\mathbb{P}\left[\mathcal{N}_{\mu, \vartheta }^{\zeta ,\theta}(t)=n\right]\\
		&=\sum_{j=k}^{\infty}\frac{(\zeta )_{j}\mathit{\Gamma } ( \vartheta )(\lambda_{\theta}G_Y(x)t)^j}{j!}E_{\mu, \vartheta+\mu j }^{\zeta+j }(-\lambda_{\theta}G_Y(x)t^{\theta})\\
		&= \mathbb{P}\left[\tilde{\mathcal{N}}_{\mu, \vartheta }^{\zeta ,\theta}(t)\geq k\right].
	\end{align*}
	Thus, we obtain the desired result.
\end{proof}

\subsection{Convergence results}
We now discuss the asymptotic behavior of the $\mathcal{N}_{\mu, \vartheta }^{\zeta ,\theta}(t),\; t\geq 0$. % as the parameters increase, we now concentrate on a feature associated with their asymptotic behavior.
\begin{proposition}\label{pooo}
	Let $\mu \in (0,1]$. Then for a fixed $0<\theta\leq 1$ and $t>0$ we have
	\begin{equation*}
		\frac{\mathcal{N}_{\mu, \vartheta }^{\zeta ,\theta}(t)}{\mathbb{E}\left[\mathcal{N}_{\mu, \vartheta }^{\zeta ,\theta}(t)\right]} \xrightarrow[\text{prob}]{\lambda_{\theta}\to \infty} 1.
	\end{equation*}
\end{proposition}
\begin{proof}
	First, we consider the convergence in mean of the random variable $ \frac{\mathcal{N}_{\mu, \vartheta }^{\zeta ,\theta}(t)}{\mathbb{E}\left[\mathcal{N}_{\mu, \vartheta }^{\zeta ,\theta}(t)\right]}$ to 1. By virtue of  triangular inequality, we obtain
	$$
	\mathbb{E}\left[\left|\frac{\mathcal{N}_{\mu, \vartheta }^{\zeta ,\theta}(t)}{\mathbb{E}\left[\mathcal{N}_{\mu, \vartheta }^{\zeta ,\theta}(t)\right]}-1\right|\right] \leq 2.
	$$
	By applying the dominated convergence theorem and evaluating the following limit
	\begin{equation}\label{conv}
		\lim _{\lambda_{\theta} \rightarrow\infty} \mathbb{E}\left[\left| \frac{\mathcal{N}_{\mu, \vartheta }^{\zeta ,\theta}(t)}{\mathbb{E}\left[\mathcal{N}_{\mu, \vartheta }^{\zeta ,\theta}(t)\right]}-1\right|\right]=\lim _{\lambda_{\theta} \rightarrow\infty} \sum_{j=0}^{\infty}\left|\frac{j}{\frac{\zeta  \mathit{\Gamma } ( \vartheta )\lambda_{\theta}t^{\theta}}{\mathit{\Gamma }(\mu+ \vartheta )}}-1\right| \frac{\left(\lambda_{\theta} t^{\theta}\right)^j\mathit{\Gamma } ( \vartheta )(\zeta )_{j}}{j!} E_{\mu, \vartheta +j\mu}^{\zeta +j}(-\lambda_{\theta}t^{\theta}).
	\end{equation}
	Taking into consideration the behavior of $E_{\alpha,\iota}^{\xi }(z)$ for large $z$ is  (see \cite{saxena2004unified}) 
	$$
	E_{\alpha,\iota}^{\xi }(z)\sim \mathscr{O}\left(|z|^{-\xi}\right), \quad|z|>1.
	$$
	This implies that the limit in Eq. (\ref{conv}) equals 0. Therefore, the proposition is proved as convergence in probability is implied by convergence in mean.
\end{proof}

\begin{proposition}\label{poi}
	Let $\mu \in (0,1]$ and $p\in\mathbb{N}$. Then, for a fixed $0<\theta\leq 1$ and $t>0$ we have
	\begin{equation*}
		\frac{\left[\mathcal{N}_{\mu, \vartheta }^{\zeta ,\theta}(t)\right]^p}{\mathbb{E}\left[\mathcal{N}_{\mu, \vartheta }^{\zeta ,\theta}(t)\right]^p} \xrightarrow[\text{prob}]{\lambda_{\theta}\to \infty} 1.
	\end{equation*}
\end{proposition}
\begin{proof}
	Using the expression given in Eq. 98 of \cite{laskin2024new} i.e.
	\begin{equation*}
		\mathbb{E}\left[\mathcal{N}_{\mu, \vartheta }^{\zeta ,\theta}(t)\right]^p= \sum_{i=0}^{p}S_{\mu, \vartheta }^{\zeta }(p,i)(\lambda_{\theta}t^{\theta})^p,
	\end{equation*}
	where $S_{\mu, \vartheta }^{\zeta }(p,i)$ is the generalized fractional Stirling number of the second kind defined in \cite{laskin2024new}. The proof of this Proposition follows on the similar lines as of Proposition \ref{pooo}. 
\end{proof}
\begin{remark}
	If we set $\theta=\mu$  and $ \vartheta =\zeta =1$ in Propositions \ref{pooo} and \ref{poi}, our results correspond to the convergence law of the FPP (see \cite{di2015fractional}).
\end{remark}
\section{Time-changed fractional counting process}\label{sec:tcfcp}
In this section, we present various distributional properties of the FCP time changed by L\'{e}vy subordinator. Next, we define its multiplicative and compound additive variants along with some of their distributional properties and typical examples.
\subsection{Time-changed fractional counting process and its distributional properties}
\begin{definition} Let  $\mathcal{N}_{\mu, \vartheta }^{\zeta ,\theta}(t),\; t\geq 0$ be the FCP defined in Eq. (\ref{process}). The time-changed fractional counting process (TCFCP), denoted by $Z(t,\lambda _{\theta }),\;t\geq 0$  is obtained by the subordination of a non-decreasing L\'{e}vy process to an independent FCP and is given by
\begin{equation*}
    Z(t,\lambda _{\theta })= N_{\mu, \vartheta }^{\zeta ,\theta}(H(t)),\;\;\; t\geq 0.
\end{equation*}
\end{definition}
\noindent The pmf of $Z(t,\lambda _{\theta }),\;t\geq 0$, denoted by $z(n,t)=\mathbb{P}\{Z(t,\lambda _{\theta })=n\}$ is derived as
\begin{align*}
    z(n,t)&= \int_{0}^{\infty}\mathcal{P}_{\mu, \vartheta }^{\zeta ,\theta}(n,y)h(y,t)\mathrm{d} y\\
    &=  (\zeta )_n \mathit{\Gamma } ( \vartheta )\frac{\lambda_{\theta}^n}{n!}\sum_{k=0}^{\infty}\frac{(-\lambda_{\theta})^k(\zeta +n)_{k}}{k!\mathit{\Gamma } (\mu(n+k)+ \vartheta )}\int_{0}^{\infty} y^{\theta(n+k)}h(y,t)\mathrm{d} y.
\end{align*}
Hence, the pmf is given by
\begin{align}\label{pmfnew}
     z(n,t)&= (\zeta )_n \mathit{\Gamma } ( \vartheta )\frac{\lambda_{\theta}^n}{n!}\sum_{k=0}^{\infty}\frac{(-\lambda_{\theta})^k(\zeta +n)_{k}}{k!\mathit{\Gamma } (\mu(n+k)+ \vartheta )}\mathbb{E}\left[(H(t))^{\theta(n+k)}\right], 
\end{align}
where $\mathbb{E}[(H(t))^{\alpha}]$ denotes the $\alpha>0$ order moment of L\'{e}vy subordinator.  A power series for the probability distribution  $z(n,x)$ is given by
\begin{align}\label{pd}
    z(n,x)&= \frac{(\zeta )_n \mathit{\Gamma } ( \vartheta )}{n!}\sum_{k=0}^{\infty}\frac{(-1)^k(\zeta +n)_{k}}{k!\mathit{\Gamma } (\mu(n+k)+ \vartheta )}\mathbb{E}\left[X^{\theta(n+k)}\right],
\end{align}
where $X$ denotes a random variable whose all fractional moments exist.

\begin{proposition}\label{proplt}
    The LT of the TCFCP $Z(t,\lambda _{\theta }),\;t\geq 0$ is denoted by $\mathcal{L}(s,t)$, is given as
    \begin{equation}\label{lt2}
        \mathcal{L}(s,t) = \mathit{\Gamma } ( \vartheta )\sum_{m=0}^{\infty}\frac{(\zeta )_m(\lambda_{\theta}(e^{-s}-1))^m}{m!\mathit{\Gamma } (\mu m+ \vartheta )}\mathbb{E}\left[(H(t))^{m\theta}\right].
    \end{equation}
\end{proposition}
\begin{proof}
    Using the conditioning argument and with the help of  Eq. (\ref{LT}),  we get
\begin{align*}
    \mathbb{E}\left[e^{-sZ(t,\lambda _{\theta })}\right]&=\mathbb{E} \left [\mathbb{E}\left [e^{-sN_{\mu, \vartheta }^{\zeta ,\theta}(H(t))}\;\big{|}\;H(t)\right]\right]\\
    &= \mathbb{E} \left [\mathbb{E}\left [\mathit{\Gamma } ( \vartheta )\sum_{m=0}^{\infty}\frac{(\zeta )_m\left(\lambda_{\theta}(H(t))^{\theta}(e^{-s}-1)\right)^m}{m!\mathit{\Gamma } (\mu m+ \vartheta )}\;\big{|}\;H(t)\right]\right]\\
    &= \mathit{\Gamma } ( \vartheta )\sum_{m=0}^{\infty}\frac{(\zeta )_m(\lambda_{\theta}(e^{-s}-1))^m}{m!\mathit{\Gamma } (\mu m+ \vartheta )}\mathbb{E}\left[(H(t))^{m\theta}\right].
    \end{align*}
As a result, the proof is complete.
\end{proof}
\begin{remark}
    The pgf of the TCFCP $Z(t,\lambda _{\theta }),\;t\geq 0$, denoted by $\mathcal{G}(s,t)$, can be obtained with the help of its LT Eq. (\ref{lt2}) and is given as
    \begin{equation*}
        \mathcal{G}(s,t) = \mathit{\Gamma } ( \vartheta )\sum_{m=0}^{\infty}\frac{(\zeta )_m (\lambda_{\theta}(s-1))^m}{m!\mathit{\Gamma } (\mu m+ \vartheta )}\mathbb{E}\left[(H(t))^{m\theta}\right].
    \end{equation*}
\end{remark}
\begin{remark}
    The TCFCP $Z(t,\lambda _{\theta }),\;t\geq 0$ doesn't have independent increments property.
 This can be observed with the help of LT of the process as given in Eq. (\ref{lt2})
    \begin{align*}
        \mathcal{L}(s,t_1+t_2) &= \mathit{\Gamma } ( \vartheta )\sum_{m=0}^{\infty}\frac{(\zeta )_m(\lambda_{\theta}(e^{-s}-1))^m}{m!\mathit{\Gamma } (\mu m+ \vartheta )}\mathbb{E}\left[(H(t_1+t_2))^{m\theta}\right]\\
        & = \mathit{\Gamma } ( \vartheta )\sum_{m=0}^{\infty}\frac{(\zeta )_m (\lambda_{\theta}(e^{-s}-1))^m}{m!\mathit{\Gamma } (\mu m+ \vartheta )}\mathbb{E}\left[(H(t_1))^{m\theta}\right]\mathbb{E}\left[(H(t_2))^{m\theta}\right].
        \end{align*}
       From this we get
        \begin{align*}
            \mathcal{L}(s,t_1+t_2) \neq \mathcal{L}(s,t_1) \mathcal{L}(s,t_2).
        \end{align*}
        It can be observed that although a Lévy subordinator exhibits independent increments, the TCFCP does not possess independent increments.  
\end{remark}
\begin{proposition}
    The mgf of the TCFCP $Z(t,\lambda _{\theta }),\;t\geq 0$,
    denoted by $\mathcal{H}(s,t)$, is given as
    \begin{equation*}
        \mathcal{H}(s,t)= \mathit{\Gamma } ( \vartheta )\sum_{m=0}^{\infty}\frac{(\zeta )_m (\lambda_{\theta}(e^{-s}-1))^m}{m!\mathit{\Gamma } (\mu m+ \vartheta )}\mathbb{E}\left[(H(t))^{m\theta}\right].
    \end{equation*}
    \end{proposition}
    \begin{proof}
        Note that
        \begin{equation*}
            \mathcal{H}(s,t)= \sum_{n=0}^{\infty}e^{-sn}z(n,t).        \end{equation*}
            Hence, the moment of any integer order can be found as
            \begin{equation*}
                h^k(s,t)=(-1)^k \frac{\partial^k}{\partial x^k}\mathcal{H}(s,t)\big{|}_{s=0}=\sum_{n=0}^{\infty}n^kz(n,t).
            \end{equation*}
       With the help of Eq. (\ref{mgf}) and the conditioning argument and by imitating the procedure as done in Proposition \ref{proplt} we obtain 
       \begin{equation*}
           \sum_{n=0}^{\infty}e^{-sn}z(n,t)=\mathit{\Gamma } ( \vartheta )\sum_{m=0}^{\infty}\frac{(\zeta )_m (\lambda_{\theta}(e^{-s}-1))^m}{m!\mathit{\Gamma } (\mu m+ \vartheta )}\mathbb{E}\left[(H(t))^{m\theta}\right].
       \end{equation*}
       Therefore, proved.\qedhere
       \end{proof}
       By replacing $\lambda_{\theta}(H(t))^{\theta}$ by the random variable $X$ as done in \cite{laskin2024new}, we obtain the mgf for the probability distribution defined in Eq. (\ref{pd}). We present this expression for later use.
       \begin{equation}\label{mgsx}
            \mathcal{H}_1(s,x)= \mathit{\Gamma } ( \vartheta )\sum_{m=0}^{\infty}\frac{(\zeta )_m}{m!\mathit{\Gamma } (\mu m+ \vartheta )}\mathbb{E}\left[(X^{\theta}(e^{-s}-1))^m\right].\qedhere
        \end{equation}
    
    \begin{theorem}
    	The mean and variance of the TCFCP $Z(t,\lambda _{\theta }),\;t\geq 0$ are\\
    	\noindent (i) $\mathbb{E}\left [ Z(t,\lambda _{\theta }) \right ]= \frac{\zeta  \mathit{\Gamma } ( \vartheta )\lambda_{\theta}}{\mathit{\Gamma }(\mu+ \vartheta )}\mathbb{E}\left[(H(t))^{\theta}\right]$\\
    	\noindent (ii) $\text{Var}\left [ Z(t,\lambda _{\theta }) \right ]= \mathbb{E}\left [ Z(t,\lambda _{\theta }) \right ]\left [ 1- \mathbb{E}\left [ Z(t,\lambda _{\theta }) \right ] \right ] + \left [\left (  \frac{ \zeta  \mathit{\Gamma } ( \vartheta )\lambda_{\theta}}{\mathit{\Gamma }(\mu+ \vartheta )}\right)^2 \left(1+\frac{1}{\zeta }\right)\frac{B(\mu+ \vartheta ,\mu+ \vartheta )}{B(2\mu+ \vartheta , \vartheta )}\right ]\mathbb{E}\left[(H(t))^{2\theta}\right].$
    \end{theorem}
    \begin{proof} With the help of Eq. (\ref{mean}) and Eq. (\ref{var}) and conditioning arguments, we have
    	\begin{align*}
    		\mathbb{E}\left [  Z(t,\lambda _{\theta }) \right ]&=\mathbb{E}\left [N_{\mu, \vartheta }^{\zeta ,\theta}(H(t))\right ]= \mathbb{E}\left [ N_{\mu, \vartheta }^{\zeta ,\theta}(H(t))\;\big{|}\;H(t) \right ]\\
    		&=\mathbb{E}\left [ \zeta  \mathit{\Gamma } ( \vartheta )\frac{\lambda_{\theta}{(H(t))}^{\theta}}{\mathit{\Gamma }(\mu+ \vartheta )}\right ] = \frac{ \zeta  \mathit{\Gamma } ( \vartheta )\lambda_{\theta}}{\mathit{\Gamma }(\mu+ \vartheta )}\mathbb{E}\left[(H(t))^{\theta}\right].
    	\end{align*}
    	The variance of the  TCFCP $ Z(t,\lambda _{\theta }),\;t\geq 0$ is evaluated as
    	\begin{align*}
    		\text{Var}\left [ Z(t,\lambda _{\theta }) \right ]&=\text{Var}\left[ \mathbb{E}\left [N_{\mu, \vartheta }^{\zeta ,\theta}(H(t))\;\big{|}\;H(t)\right ] \right] +\mathbb{E}\left [\text{Var}\left[N_{\mu, \vartheta }^{\zeta ,\theta}(H(t))\;\big{|}\;H(t) \right ]\right]\\
    		&=\text{Var}\left[\frac{\zeta  \mathit{\Gamma } ( \vartheta )\lambda_{\theta}{(H(t))}^{\theta}}{\mathit{\Gamma }(\mu+ \vartheta )}\right] \\
    		& \;\;+ \mathbb{E}\left [\frac{ \zeta  \mathit{\Gamma } ( \vartheta )\lambda_{\theta}}{\mathit{\Gamma }(\mu+ \vartheta )}{(H(t))}^{\theta}+\left (  \frac{ \zeta  \mathit{\Gamma } ( \vartheta )\lambda_{\theta}}{\mathit{\Gamma }(\mu+ \vartheta )}\right)^2 \left[\left(1+\frac{1}{\zeta }\right)\frac{B(\mu+ \vartheta ,\mu+ \vartheta )}{B(2\mu+ \vartheta , \vartheta )}-1\right ]{(H(t))}^{2\theta}\right]\\
    		&=  \left (\frac{ \zeta  \mathit{\Gamma } ( \vartheta )\lambda_{\theta}}{\mathit{\Gamma }(\mu+ \vartheta )}  \right )^2\left [\mathbb{E}\left[(H(t))^{2\theta}\right] - \left[\mathbb{E}\left[(H(t))^{\theta}\right] \right ]^2\right] + \frac{ \zeta  \mathit{\Gamma } ( \vartheta )\lambda_{\theta}}{\mathit{\Gamma }(\mu+ \vartheta )}\mathbb{E}\left[(H(t))^{\theta}\right]\\
    		& \;\;+\left (  \frac{ \zeta  \mathit{\Gamma } ( \vartheta )\lambda_{\theta}}{\mathit{\Gamma }(\mu+ \vartheta )}\right)^2 \left[\left(1+\frac{1}{\zeta }\right)\frac{B(\mu+ \vartheta ,\mu+ \vartheta )}{B(2\mu+ \vartheta , \vartheta )}-1\right ]\mathbb{E}\left[(H(t))^{2\theta}\right].\qedhere
    	\end{align*}
    \end{proof}

\begin{theorem}
    The waiting time distribution of the TCFCP $Z(t,\lambda _{\theta }),\;t\geq 0$ is given as
    \begin{align*}
    \phi (\tau ) &= -\frac{\mathrm{d} }{\mathrm{d} \tau }\left [\mathit{\Gamma } ( \vartheta )\sum_{k=0}^{\infty}\frac{(-\lambda_{\theta})^k(\zeta )_{k}}{k!\mathit{\Gamma } (\mu k+ \vartheta )}\mathbb{E}\left[(H(\tau))^{\theta k}\right]  \right ],\tau>0.
\end{align*}
\end{theorem}
 \begin{proof} Let us denote the waiting time distribution  by $\phi (\tau ),\tau>0$. Consider
\begin{equation}\label{wt1}
    \phi (\tau ) = -\frac{\mathrm{d} }{\mathrm{d} \tau }z(\tau),
\end{equation}
where the probability of a given waiting time being greater or equal to $\tau$ is given by 
\begin{equation}\label{wt2}
    z(\tau) = 1- \sum_{n=1}^{\infty}z(n,\tau) = z(0,\tau) = \mathit{\Gamma } ( \vartheta )\sum_{k=0}^{\infty}\frac{(-\lambda_{\theta})^k(\zeta )_{k}}{k!\mathit{\Gamma } (\mu k+ \vartheta )}\mathbb{E}\left[(H(\tau))^{\theta k}\right].
\end{equation}
Using Eqs. (\ref{wt1}) and (\ref{wt2}), we evaluate the waiting time distribution as
\begin{align*}
    \phi (\tau ) &= -\frac{\mathrm{d} }{\mathrm{d} \tau }\left [\mathit{\Gamma } ( \vartheta )\sum_{k=0}^{\infty}\frac{(-\lambda_{\theta})^k(\zeta )_{k}}{k!\mathit{\Gamma } (\mu k+ \vartheta )}\mathbb{E}\left[(H(\tau))^{\theta k}\right]  \right ].
\end{align*}
Consequently, the theorem has been proved.
\end{proof}
\begin{theorem}
The first passage time distribution of the TCFCP $Z(t,\lambda _{\theta }),\;t\geq 0$ denoted by $\mathcal{T}_w$ for $t\geq0,\; w\in \mathbb{N}$ is given as
    \begin{equation*}
    	\mathbb{P}\left[\mathcal{T}_w>t\right] = \sum_{n=0}^{w-1}(\zeta )_n \mathit{\Gamma } ( \vartheta )\frac{\lambda_{\theta}^n}{n!}\sum_{k=0}^{\infty}\frac{(-\lambda_{\theta})^k(\zeta +n)_{k}}{k!\mathit{\Gamma } (\mu(n+k)+ \vartheta )}\mathbb{E}\left[(H(t))^{\theta(n+k)}\right].
    \end{equation*}
    \end{theorem}
\begin{proof}
The first passage time is the time at which, for the first time, the process reaches a certain threshold. Let $\mathcal{T}_w$ be the time of first upcrossing of level $w$ and is defined as
$	\mathcal{T}_w := \inf\{t\geq0 :z(n,t)\geq w \}.$ Therefore, the survival function $\mathbb{P}\left[\mathcal{T}_w>t\right]$ can be derived as
\begin{align*}
\mathbb{P}\left[\mathcal{T}_w>t\right] &= \mathbb{P}\left[z(n,t)<w \right]
	= \sum_{n=0}^{w-1}\mathbb{P}\left[z(n,t) \right]\\
	&=\sum_{n=0}^{w-1}\int_{0}^{\infty}\mathcal{P}_{\mu, \vartheta }^{\zeta ,\theta}(n,y)h(y,t)\mathrm{d} y\\
	&=\sum_{n=0}^{w-1}(\zeta )_n \mathit{\Gamma } ( \vartheta )\frac{\lambda_{\theta}^n}{n!}\sum_{k=0}^{\infty}\frac{(-\lambda_{\theta})^k(\zeta +n)_{k}}{k!\mathit{\Gamma } (\mu(n+k)+ \vartheta )}\int_{0}^{\infty} y^{\theta(n+k)}h(y,t)\mathrm{d} y\\
	&=\sum_{n=0}^{w-1}(\zeta )_n \mathit{\Gamma } ( \vartheta )\frac{\lambda_{\theta}^n}{n!}\sum_{k=0}^{\infty}\frac{(-\lambda_{\theta})^k(\zeta +n)_{k}}{k!\mathit{\Gamma } (\mu(n+k)+ \vartheta )}\mathbb{E}\left[(H(t))^{\theta(n+k)}\right].
\end{align*}
Moreover, the distribution of $\mathcal{T}_w$ can be written as
\begin{align*}
    \mathbb{P}\left[\mathcal{T}_w<t\right] &= \mathbb{P}\left[z(n,t)\geq w \right]
	= \sum_{n=w}^{\infty}\mathbb{P}\left[z(n,t) \right]\\
 &=\sum_{n=w}^{\infty}(\zeta )_n \mathit{\Gamma } ( \vartheta )\frac{\lambda_{\theta}^n}{n!}\sum_{k=0}^{\infty}\frac{(-\lambda_{\theta})^k(\zeta +n)_{k}}{k!\mathit{\Gamma } (\mu(n+k)+ \vartheta )}\mathbb{E}\left[(H(t))^{\theta(n+k)}\right].
\end{align*} 
Therefore, the density function $\mathcal{P}(n,t)= \mathbb{P}\left[\mathcal{T}_k \in dt \right]/dt$ is
\begin{align*}
	\mathcal{P}(n,t)&= \frac{d}{dt}\sum_{n=w}^{\infty}(\zeta )_n \mathit{\Gamma } ( \vartheta )\frac{\lambda_{\theta}^n}{n!}\sum_{k=0}^{\infty}\frac{(-\lambda_{\theta})^k(\zeta +n)_{k}}{k!\mathit{\Gamma } (\mu(n+k)+ \vartheta )}\mathbb{E}\left[(H(t))^{\theta(n+k)}\right]\\
	&=\frac{d}{dt}\left ( 1-\sum_{n=0}^{w-1}(\zeta )_n \mathit{\Gamma } ( \vartheta )\frac{\lambda_{\theta}^n}{n!}\sum_{k=0}^{\infty}\frac{(-\lambda_{\theta})^k(\zeta +n)_{k}}{k!\mathit{\Gamma } (\mu(n+k)+ \vartheta )}\mathbb{E}\left[(H(t))^{\theta(n+k)}\right]\right )\\
	&= -\frac{d}{dt}\sum_{n=0}^{w-1}(\zeta )_n \mathit{\Gamma } ( \vartheta )\frac{\lambda_{\theta}^n}{n!}\sum_{k=0}^{\infty}\frac{(-\lambda_{\theta})^k(\zeta +n)_{k}}{k!\mathit{\Gamma } (\mu(n+k)+ \vartheta )}\mathbb{E}\left[(H(t))^{\theta(n+k)}\right].\qedhere
\end{align*}
\end{proof}

\subsection{Multiplicative compound fractional counting process} 
In this subsection, we study the multiplicative compound fractional counting process (MCFCP) denoted by $\{\mathcal{N}_{\pi}(t),\;t\geq 0\}$ where the 
randomized time is governed by the FCP. It is given by
\begin{equation*}
    \mathcal{N}_{\pi}(t) = \prod_{j=1}^{\mathcal{N}_{\mu, \vartheta }^{\zeta ,\theta}(t)}X_j,
\end{equation*}
where $\{X_n,\; n\geq 1\}$ is a sequence of random variables independent of the FCP. Let us set
\begin{equation*}
    R_m = X_1X_2\cdots X_m.
\end{equation*}
The cdf $W_{\mathcal{N}_{\pi}}(y,t)$ of $\mathcal{N}_{\pi}(t)$, for $y\in \mathbb{R}$ and $t>0$, can be evaluated as
\begin{align*}
    W_{\mathcal{N}_{\pi}}(y,t)& = \sum_{m=0}^{\infty}\mathbb{P}\left[\mathcal{N}_{\mu, \vartheta }^{\zeta ,\theta}(t)=m\right]\mathbb{P}\left[R_m\leq y\right]\\
    &= \mathbb{I}_{\{ y \geq 1\}}\mathcal{P}_{\mu, \vartheta }^{\zeta ,\theta}(0,t) +  \sum_{m=1}^{\infty}\mathbb{P}\left[\mathcal{N}_{\mu, \vartheta }^{\zeta ,\theta}(t)=m\right]\mathbb{P}\left[R_m\leq y\right],
\end{align*}
with the help of Eq. (\ref{process}), we get the following 
\begin{equation}\label{cdf1}
     W_{\mathcal{N}_{\pi}}(y,t)= \mathbb{I}_{\{ y \geq 1\}} \mathit{\Gamma } ( \vartheta )E_{\mu, \vartheta }^{\zeta }(-\lambda_{\theta}t^{\theta})+  \sum_{m=1}^{\infty}\mathbb{P}\left[\mathcal{N}_{\mu, \vartheta }^{\zeta ,\theta}(t)=m\right]\mathbb{P}\left[R_m\leq y\right].
\end{equation}
If $R_m$ is absolutely continuous having probability density $f_{R_m}(\cdot)$ then due to Eq. (\ref{cdf1}), then we have the following additional observations\\
(i) The probability law of $\mathcal{N}_{\pi}$ has an absolutely continuous component that is expressed by the density as
\begin{equation*}
    h_{\mathcal{N}_{\pi}}(y,t) = \sum_{m=1}^{\infty} \mathbb{P}\left[\mathcal{N}_{\mu, \vartheta }^{\zeta ,\theta}(t)=m\right]f_{R_m}(y);\;\;\; y\neq 1,\; t>0.
\end{equation*}
(ii) The discrete component is given by 
\begin{align*}
    \mathbb{P}\left[\mathcal{N}_{\pi}(t)=1\right] = \mathbb{P}\left[\mathcal{N}_{\mu, \vartheta }^{\zeta ,\theta}(t)= 0\right] &=  W_{\mathcal{N}_{\pi}}(1,t) -  W_{\mathcal{N}_{\pi}}(y,t),\; y<1\\
    &=  \mathit{\Gamma } ( \vartheta )E_{\mu, \vartheta }^{\zeta }(-\lambda_{\theta}t^{\theta}).
\end{align*}
\begin{remark}
When $\mu= \vartheta =\zeta =\theta=1$, the MCFCP reduces to the multiplicative compound Poisson process considered in \cite{meoli2023some}.
\end{remark}
\noindent 
Now we consider the MCFCP at L\'{e}vy times denoted by $Z_{\pi} = \{Z_{\pi}(t),\; t\geq 0\}$ defined as 
\begin{equation*}
    Z_{\pi}(t) = \{\mathcal{N}_{\pi}(H(t))\}, t\geq 0,
\end{equation*}
where $H(t)$ is an independent L\'{e}vy subordinator. From the definition, it follows that if $\mathcal{N}_{\mu, \vartheta }^{\zeta ,\theta}(t)= 0$ then $\mathcal{N}_{\pi}(t)=1$, therefore we can say that $\mathbb{P}\left[\mathcal{N}_{\pi}(t)=1\right] \geq \mathbb{P}\left[\mathcal{N}_{\mu, \vartheta }^{\zeta ,\theta}(t)= 0\right] $.

We next evaluate the probability law of $Z_{\pi}$.
\begin{proposition}
    Let $y\in \mathbb{R}$ and $t>0$, the cdf $ W_{Z_{\pi}}(y,t)$ of the MCFCP subordinated with an independent L\'{e}vy subordinator is 
    \begin{align*}
        W_{Z_{\pi}}(y,t) =& \;\mathbb{I}_{\{ y \geq 1\}} \mathit{\Gamma } ( \vartheta )\sum_{k=0}^{\infty}\frac{(\zeta )_k(-\lambda_{\theta})^k}{k!\mathit{\Gamma } (\mu k+ \vartheta )} \mathbb{E}\left[(H(t))^{\theta k}\right]\\ & + \mathit{\Gamma } ( \vartheta )\sum_{m=1}^{\infty}\frac{\lambda_{\theta}^m(\zeta )_m}{m!}\mathbb{P}\left[R_m\leq y\right]\sum_{l=0}^{\infty}\frac{(\zeta  + m)_{l}(-\lambda_{\theta})^l}{l! \mathit{\Gamma } (\mu (l+m)+ \vartheta )}\mathbb{E}\left[(H(t))^{\theta(l+m) }\right].\nonumber
    \end{align*}
\end{proposition}
\begin{proof} Let $y\in \mathbb{R}$ and $t>0$, we have that
    \begin{align*}
         W_{Z_{\pi}}(y,t) &= \sum_{n=0}^{\infty}\mathbb{P}\left[ H(t) = n\right]W_{\mathcal{N}_{\pi}}(y,n)\\
         &= \sum_{n=0}^{\infty}\mathbb{P}\left[ H(t) = n\right]\left [  \mathbb{I}_{\{ y \geq 1\}} \mathit{\Gamma } ( \vartheta )E_{\mu, \vartheta }^{\zeta }(-\lambda_{\theta}n^{\theta})+  \sum_{m=1}^{\infty}\mathbb{P}\left[\mathcal{N}_{\mu, \vartheta }^{\zeta ,\theta}(n)=m\right]\mathbb{P}\left[R_m\leq y\right]\right ].
    \end{align*}
    Evaluating the above expression in two parts. Consider the first part as follows
    \begin{align*}
        \mathbb{I}_{\{ y \geq 1\}} \mathit{\Gamma } ( \vartheta )\sum_{n=0}^{\infty}&\mathbb{P}\left[H(t) = n\right]E_{\mu, \vartheta }^{\zeta }(-\lambda_{\theta}n^{\theta}) \\
        &=\mathbb{I}_{\{ y \geq 1\}} \mathit{\Gamma } ( \vartheta )\sum_{n=0}^{\infty}\mathbb{P}\left[ H(t) = n\right]\sum_{k=0}^{\infty}\frac{(\zeta )_k(-\lambda_{\theta}n^{\theta})^k}{k!\mathit{\Gamma } (\mu k+ \vartheta )}\\
        &=\mathbb{I}_{\{ y \geq 1\}} \mathit{\Gamma } ( \vartheta )\sum_{k=0}^{\infty}\frac{(\zeta )_k(-\lambda_{\theta})^k}{k!\mathit{\Gamma } (\mu k+ \vartheta )}\sum_{n=0}^{\infty} n^{\theta k} \mathbb{P}\left[ H(t) = n\right]\\
         &=\mathbb{I}_{\{ y \geq 1\}} \mathit{\Gamma } ( \vartheta )\sum_{k=0}^{\infty}\frac{(\zeta )_k(-\lambda_{\theta})^k}{k!\mathit{\Gamma } (\mu k+ \vartheta )} \mathbb{E}\left[(H(t))^{\theta k}\right].
    \end{align*}
Now, consider the second part
\begin{align*}
    \sum_{n=0}^{\infty}&\mathbb{P}\left[ H(t) = n\right]\sum_{m=1}^{\infty}\mathbb{P}\left[\mathcal{N}_{\mu, \vartheta }^{\zeta ,\theta}(n)=m\right]\mathbb{P}\left[R_m\leq y\right]\\
    &=\mathit{\Gamma } ( \vartheta )\sum_{n=0}^{\infty}\mathbb{P}\left[ H(t) = n\right]\sum_{m=1}^{\infty}\frac{(\lambda_{\theta}n^{\theta})^m(\zeta )_m}{m!} E_{\mu, \vartheta +m\mu}^{\zeta +m}(-\lambda_{\theta}n^{\theta})\mathbb{P}\left[R_m\leq y\right]\\
    &= \mathit{\Gamma } ( \vartheta )\sum_{m=1}^{\infty}\frac{\lambda_{\theta}^m(\zeta )_m}{m!}\mathbb{P}\left[R_m\leq y\right]\sum_{n=0}^{\infty} n^{\theta m}\sum_{l=0}^{\infty}\frac{(\zeta  + m)_{l}(-\lambda_{\theta}n^\theta)^l}{l! \mathit{\Gamma } (\mu (l+m)+ \vartheta )}\mathbb{P}\left[ H(t) = n\right]\\
     &= \mathit{\Gamma } ( \vartheta )\sum_{m=1}^{\infty}\frac{\lambda_{\theta}^m(\zeta )_m}{m!}\mathbb{P}\left[R_m\leq y\right]\sum_{l=0}^{\infty}\frac{(\zeta  + m)_{l}(-\lambda_{\theta})^l}{l! \mathit{\Gamma } (\mu (l+m)+ \vartheta )}\sum_{n=0}^{\infty} n^{\theta (m+l)}\mathbb{P}\left[H(t) = n\right]\\
      &= \mathit{\Gamma } ( \vartheta )\sum_{m=1}^{\infty}\frac{\lambda_{\theta}^m(\zeta )_m}{m!}\mathbb{P}\left[R_m\leq y\right]\sum_{l=0}^{\infty}\frac{(\zeta  + m)_{l}(-\lambda_{\theta})^l}{l! \mathit{\Gamma } (\mu (l+m)+ \vartheta )}\mathbb{E}\left[(H(t))^{\theta(l+m) }\right].
\end{align*}
With the help of above two parts, the proof is complete.
\end{proof}

\begin{corollary}
If $R_m$ is absolutely continuous with density $f_{R_m}(\cdot)$,   then \\
(i) the density of $Z_{\pi}$ for any $t>0$ and $y\neq 0$ is
\begin{equation}\label{hss}
    h_{Z_{\pi}}(y,t) =  \mathit{\Gamma } ( \vartheta )\sum_{m=1}^{\infty}\frac{\lambda_{\theta}^m(\zeta )_m}{m!}\sum_{l=0}^{\infty}\frac{(\zeta  + m)_{l}(-\lambda_{\theta})^l}{l! \mathit{\Gamma } (\mu (l+m)+ \vartheta )}\mathbb{E}\left[(H(t))^{\theta(l+m) }\right]f_{R_m}(y),
\end{equation}
(ii) the probability law of $Z_{\pi}$ has a discrete component given by 
\begin{align*}
    \mathbb{P}\left[Z_{\pi}(t)=1\right] &= \sum_{n=0}^{\infty}\mathbb{P}\left[\mathcal{N}_{\pi}(n)= 1\right]\mathbb{P}\left[ H(t) = n\right]\\
    &=  \sum_{n=0}^{\infty}\mathit{\Gamma } ( \vartheta )E_{\mu, \vartheta }^{\zeta }(-\lambda_{\theta}n^{\theta})\mathbb{P}\left[ H(t) = n\right]\\
    &= \sum_{m=0}^{\infty}\frac{(\zeta )_m\mathit{\Gamma } ( \vartheta )(-\lambda_{\theta})^m}{m!\mathit{\Gamma } (\mu m+ \vartheta )}\mathbb{E}\left[(H(t))^{\theta m}\right].
\end{align*}
\end{corollary}
\begin{corollary}\label{cor2}
    If $\{X_n, n\geq 1\}$ are discrete, we set
    \begin{equation*}
        \tilde{b}_s^{*m} = \mathbb{P}\left[ X_1X_2\cdots X_m= s\right],\;\; s\in\mathbb{N}_{0}.
    \end{equation*}
    Then for $s\in \mathbb{N}_{0}$, we have that
    \begin{align}\label{corr2}
        \mathbb{P}\left[Z_{\pi}(t)=s\right]= &\;\mathbb{I}_{\{ s= 1\}} \sum_{k=0}^{\infty}\frac{(\zeta )_k\mathit{\Gamma } ( \vartheta )(-\lambda_{\theta})^k}{k!\mathit{\Gamma } (\mu k+ \vartheta )} \mathbb{E}\left[(H(t))^{\theta k}\right]\\ & + \mathit{\Gamma } ( \vartheta )\sum_{m=1}^{\infty}\frac{\lambda_{\theta}^m(\zeta )_m}{m!}\tilde{b}_s^{*m} \sum_{l=0}^{\infty}\frac{(\zeta  + m)_{l}(-\lambda_{\theta})^l}{l! \mathit{\Gamma } (\mu (l+m)+ \vartheta )}\mathbb{E}\left[(H(t))^{\theta(l+m) }\right],\nonumber
    \end{align}
    where $\mathbb{I}_{\{\cdot\}}$ is indicator function.
\end{corollary}
Next, we present the following special cases for the MCFCP at L\'evy times. 
\begin{example}
    Let $\{X_n, n\geq 1\}$ be iid Bernoulli random variables having parameter $p\in[0,1]$. Then the $R_m$ follows Bernoulli having parameter $p^m$. Therefore, in accordance to Corollary \ref{cor2}, it reduces to
    \begin{equation*}
         \tilde{b}_1^{*m} = p^m \;\;\;\text{and}\;\;\;\tilde{b}_0^{*m} = 1- p^m.
    \end{equation*}
   The Eq. (\ref{corr2}) for $s=1$, becomes
    \begin{align*}
        \mathbb{P}\{Z_{\pi}(t)=1\}= &\;\sum_{k=0}^{\infty}\frac{(\zeta )_k\mathit{\Gamma } ( \vartheta )(-\lambda_{\theta})^k}{k!\mathit{\Gamma } (\mu k+ \vartheta )} \mathbb{E}\left[(H(t))^{\theta k}\right]\\ & + \mathit{\Gamma } ( \vartheta )\sum_{m=1}^{\infty}\frac{\lambda_{\theta}^m(\zeta )_m}{m!}p^m \sum_{l=0}^{\infty}\frac{(\zeta  + m)_{l}(-\lambda_{\theta})^l}{l! \mathit{\Gamma } (\mu (l+m)+ \vartheta )}\mathbb{E}\left[(H(t))^{\theta(l+m) }\right]. 
    \end{align*}
    A similar computation leads to 
    \begin{align*}
         \mathbb{P}\{Z_{\pi}(t)=0\}=&\;\mathit{\Gamma } ( \vartheta )\sum_{m=1}^{\infty}\frac{\lambda_{\theta}^m(\zeta )_m}{m!}(1-p^m) \sum_{l=0}^{\infty}\frac{(\zeta  + m)_{l}(-\lambda_{\theta})^l}{l! \mathit{\Gamma } (\mu (l+m)+ \vartheta )}\mathbb{E}\left[(H(t))^{\theta(l+m) }\right]\\
         =&\; 1-\sum_{k=0}^{\infty}\frac{(\zeta )_k\mathit{\Gamma } ( \vartheta )(-\lambda_{\theta})^k}{k!\mathit{\Gamma } (\mu k+ \vartheta )}\mathbb{E}\left[(H(t))^{\theta k}\right]\\
         & -  \mathit{\Gamma } ( \vartheta )\sum_{m=1}^{\infty}\frac{\lambda_{\theta}^m(\zeta )_m}{m!}p^m \sum_{l=0}^{\infty}\frac{(\zeta  + m)_{l}(-\lambda_{\theta})^l}{l! \mathit{\Gamma } (\mu (l+m)+ \vartheta )}\mathbb{E}\left[(H(t))^{\theta(l+m) }\right]. 
    \end{align*}
\end{example}
\begin{example}
     Let $\{X_n, n\geq 1\}$ be iid random variables with beta distribution $\beta(c,d),c,d>0$. Then the density $f_X(\cdot)$ is given as
     \begin{equation}\label{fss}
         f_X(x)= \frac{x^{c-1}(1-x)^{d-1}}{B(c,d)}, \;\; 0<x<1.
     \end{equation}
      Fan \cite{fan1991distribution} showed that if $X_1,X_2,X_3,\cdots,X_m$ are independent random variables having beta distributions $\beta(c,d_1),\beta(c+d_1,d_2),\cdots,\beta(c+d_1+d_2+\cdots+d_{m-1},d_m)$ respectively and $c,d_1,\cdots,d_m>0$ then $R_m= X_1X_2\cdots X_m$ has Beta distribution $\beta(c,d_1+d_2+\cdots+d_m)$. In particular, if $c=d_1=d_2=\cdots=d_m=1$, Eq. (\ref{fss}) becomes 
     \begin{equation*}
         f_{R_m}(x)= \frac{(1-x)^{m-1}}{B(1,m)}, \;\; 0<x<1.
     \end{equation*}
     By substituting the above expression in  Eq. (\ref{hss}), we get the following density function
     \begin{equation*}
    h_{Z_{\pi}}(y,t) =  \mathit{\Gamma } ( \vartheta )\sum_{m=1}^{\infty}\frac{\lambda_{\theta}^m(\zeta )_m}{m!}\sum_{l=0}^{\infty}\frac{(\zeta  + m)_{l}(-\lambda_{\theta})^l}{l! \mathit{\Gamma } (\mu (l+m)+ \vartheta )}\frac{(1-y)^{m-1}}{B(1,m)}\mathbb{E}\left[(H(t))^{\theta(l+m) }\right].
\end{equation*}
\end{example}

\subsection{Fractional generalized compound process at L\'{e}vy times}
The fractional generalized compound process (FGCP), as defined in \cite{laskin2024new}, is given as 
\begin{equation*}
    \mathcal{X}_{\mu, \vartheta }^{\zeta ,\theta}(t) =  \sum_{i=1}^{\mathcal{N}_{\mu, \vartheta }^{\zeta ,\theta}(t)}Y_i,
\end{equation*}
where $\{Y_i, i\in \mathbb{N}\}$ be iid random variables, independent of the FCP, having a common distribution function $F_Y(\cdot)$. With the help of conditioning argument the cdf $W_X(y,t)$ of $ \mathcal{X}_{\mu, \vartheta }^{\zeta ,\theta}(t) $ can be evaluated as
\begin{equation*}
    W_X(y,t) = \mathbb{P}\left[\mathcal{X}_{\mu, \vartheta }^{\zeta ,\theta}(t)\leq y\right] =  \sum_{m=0}^{\infty} \mathbb{P}\left[\mathcal{N}_{\mu, \vartheta }^{\zeta ,\theta}(t)=m\right] S_Y^{(m)}(y),\;\;\;y\in \mathbb{R},\; t>0,
\end{equation*}
where $\mathbb{P}\left[\mathcal{N}_{\mu, \vartheta }^{\zeta ,\theta}(t)=m\right]$ denotes the state probabilities of the FCP. Further, for $m\in \mathbb{N}$ 
\begin{align*}
    S_Y^{(m)}(y) &= \mathbb{P}\left[Y_1+Y_2+\cdots+Y_m\leq y\right]\\
    &= \int_{-\infty}^{\infty} S_Y^{(m-1)}(y-z) \mathrm{d} S_Y(z),
\end{align*}
is the $m$-fold convolution of $S_Y(\cdot)$ and $ S_Y^{(0)}(y) = \mathbb{I}_{\{ y \geq 0\}}$. Consider
\begin{align}
    W_X(y,t) &= \mathbb{I}_{\{ y \geq 0\}}\mathcal{P}_{\mu, \vartheta }^{\zeta ,\theta}(0,t) +  \sum_{m=1}^{\infty}\mathbb{P}\left[\mathcal{N}_{\mu, \vartheta }^{\zeta ,\theta}(t)=m\right]S_Y^{(m)}(y)\nonumber \\
    &=  \mathbb{I}_{\{ y \geq 0\}} \mathit{\Gamma } ( \vartheta )E_{\mu, \vartheta }^{\zeta }(-\lambda_{\theta}t^{\theta})+  \sum_{m=1}^{\infty}\mathbb{P}\left[\mathcal{N}_{\mu, \vartheta }^{\zeta ,\theta}(t)=m\right]S_Y^{(m)}(y).\nonumber
\end{align}
The mean of the FGCP $\mathcal{X}_{\mu, \vartheta }^{\zeta ,\theta}(t),\; t>0$ can be evaluated using conditioning argument on $\mathcal{N}_{\mu, \vartheta }^{\zeta ,\theta}(t)$ as 
\begin{align}\label{meancpp}
   \mathbb{E}\left[ \mathcal{X}_{\mu, \vartheta }^{\zeta ,\theta}(t) \right]&=\mathbb{E} \left [\mathbb{E}\left [ \mathcal{X}_{\mu, \vartheta }^{\zeta ,\theta}(t) |\mathcal{N}_{\mu, \vartheta }^{\zeta ,\theta}(t)\right]\right]\nonumber\\
   &= \mathbb{E}\left[M_{Y_i}\mathcal{N}_{\mu, \vartheta }^{\zeta ,\theta}(t)\right]\nonumber\\
   &= M_{Y_i} \frac{\zeta  \mathit{\Gamma } ( \vartheta )\lambda_{\theta}t^{\theta}}{\mathit{\Gamma }(\mu+ \vartheta )}.
\end{align}
Also, the variance of the FGCP $\mathcal{X}_{\mu, \vartheta }^{\zeta ,\theta}(t)$ is evaluated as
\begin{align}\label{varcpp}
\text{Var}\left [ \mathcal{X}_{\mu, \vartheta }^{\zeta ,\theta}(t)\right ]&=\text{Var}\left[ \mathbb{E}\left [\mathcal{X}_{\mu, \vartheta }^{\zeta ,\theta}(t)|\mathcal{N}_{\mu, \vartheta }^{\zeta ,\theta}(t)\right ] \right] +\mathbb{E}\left [\text{Var}\left[\mathcal{X}_{\mu, \vartheta }^{\zeta ,\theta}(t)|\mathcal{N}_{\mu, \vartheta }^{\zeta ,\theta}(t) \right ]\right]\nonumber\\
&=  \text{Var}\left[M_{Y_i} \mathcal{N}_{\mu, \vartheta }^{\zeta ,\theta}(t)\right] +\mathbb{E}\left[V_{Y_i}^2\mathcal{N}_{\mu, \vartheta }^{\zeta ,\theta}(t)\right]\nonumber\\
 &=(M_{Y_i}^2+V_{Y_i}^2)\mathbb{E}\left[\mathcal{N}_{\mu, \vartheta }^{\zeta ,\theta}(t)\right] +\left [M_{Y_i}\mathbb{E}\left[\mathcal{N}_{\mu, \vartheta }^{\zeta ,\theta}(t)\right] \right ]^{2}\left\{ \left(1+\frac{1}{\zeta }\right)\frac{B(\mu+ \vartheta ,\mu+ \vartheta )}{B(2\mu+ \vartheta , \vartheta )}-1\right\}\nonumber\\
 &=(M_{Y_i}^2+V_{Y_i}^2)\frac{\zeta  \mathit{\Gamma } ( \vartheta )\lambda_{\theta}t^{\theta}}{\mathit{\Gamma }(\mu+ \vartheta )} +\left [M_{Y_i}\frac{\zeta  \mathit{\Gamma } ( \vartheta )\lambda_{\theta}t^{\theta}}{\mathit{\Gamma }(\mu+ \vartheta )} \right ]^{2}\left\{ \left(1+\frac{1}{\zeta }\right)\frac{B(\mu+ \vartheta ,\mu+ \vartheta )}{B(2\mu+ \vartheta , \vartheta )}-1\right\},
\end{align} 
where $M_{Y_i}$ and $V_{Y_i}$ denotes the mean and variance of $Y_i$, respectively.

\begin{definition}[FGCP at L\'evy times]
   Let us define the fractional generalized compound process (FGCP) at L\'{e}vy times as follows
\begin{align*}
    Z(t) = \mathcal{X}_{\mu, \vartheta }^{\zeta ,\theta}(H(t)),\;\; t\geq0,
\end{align*}
where $H(t),\; t\geq 0$ denotes an independent L\'{e}vy subordinator and $\mathcal{X}_{\mu, \vartheta }^{\zeta ,\theta}(t)$ is the  FGCP. 
\end{definition}
The mean of $Z(t),\; t>0$ is evaluated with the help of conditioning argument and Eq. (\ref{meancpp})
\begin{align}
	\mathbb{E}\left[ Z(t) \right]&=\mathbb{E} \left [\mathbb{E}\left [ \mathcal{X}_{\mu, \vartheta }^{\zeta ,\theta}(H(t)) |H(t)\right]\right]\nonumber\\
	&= M_{Y_i} \frac{\zeta  \mathit{\Gamma } ( \vartheta )\lambda_{\theta}}{\mathit{\Gamma }(\mu+ \vartheta )}\mathbb{E}\left[(H(t))^{\theta}\right]\nonumber.
\end{align}
\begin{proposition}
 Let $y\in \mathbb{R}$ and $t>0$, the cdf $ H_{Z}(y,t)$ of the FGCP at L\'{e}vy times is given by
    \begin{align}
        W_{Z}(y,t) =& \;\mathbb{I}_{\{ y \geq 0\}} \sum_{k=0}^{\infty}\frac{(\zeta )_k\mathit{\Gamma } ( \vartheta )(-\lambda_{\theta})^k}{k!\mathit{\Gamma } (\mu k+ \vartheta )} \mathbb{E}\left[(H(t))^{\theta k}\right]\nonumber\\ & + \mathit{\Gamma } ( \vartheta )\sum_{m=1}^{\infty}\frac{(\lambda_{\theta})^m(\zeta )_m}{m!}S_Y^{(m)}(y)\sum_{l=0}^{\infty}\frac{(\zeta  + m)_{l}(-\lambda_{\theta})^l}{l! \mathit{\Gamma } (\mu (l+m)+ \vartheta )}\mathbb{E}\left[(H(t))^{\theta(l+m) }\right].\nonumber
    \end{align}    
\end{proposition}
\begin{proof}
    With the help of conditioning argument, we have 
\begin{align*}
         W_{Z}(y,t) &= \sum_{n=0}^{\infty}\mathbb{P}\left[H(t) = n\right]W_X(y,n)\\
         &= \sum_{n=0}^{\infty}\mathbb{P}\left[ H(t) = n\right]\left [  \mathbb{I}_{\{ y \geq 0\}} \mathit{\Gamma } ( \vartheta )E_{\mu, \vartheta }^{\zeta }(-\lambda_{\theta}n^{\theta})+  \sum_{m=1}^{\infty}\mathbb{P}\left[\mathcal{N}_{\mu, \vartheta }^{\zeta ,\theta}(n)=m\right]S_Y^{(m)}(y)\right ].
    \end{align*}
    We now evaluate the RHS of the above expression in two parts. The first part is as follows 
\begin{align*}
        \mathbb{I}_{\{ y \geq 0\}} \mathit{\Gamma } ( \vartheta )\sum_{n=0}^{\infty}&\mathbb{P}\left[ H(t) = n\right]E_{\mu, \vartheta }^{\zeta }(-\lambda_{\theta}n^{\theta}) \\
        &=\mathbb{I}_{\{ y \geq 0\}} \mathit{\Gamma } ( \vartheta )\sum_{n=0}^{\infty}\mathbb{P}\left[ H(t) = n\right]\sum_{k=0}^{\infty}\frac{(\zeta )_k(-\lambda_{\theta}n^{\theta})^k}{k!\mathit{\Gamma } (\mu k+ \vartheta )}\\
         &=\mathbb{I}_{\{ y \geq 0\}} \mathit{\Gamma } ( \vartheta )\sum_{k=0}^{\infty}\frac{(\zeta )_k(-\lambda_{\theta})^k}{k!\mathit{\Gamma } (\mu k+ \vartheta )} \mathbb{E}\left[(H(t))^{\theta k}\right].
    \end{align*}
The second part is given by
\begin{align*}
    \sum_{n=0}^{\infty}&\mathbb{P}\left[ H(t) = n\right]\sum_{m=1}^{\infty}\mathbb{P}\left[\mathcal{N}_{\mu, \vartheta }^{\zeta ,\theta}(n)=m\right]S_Y^{(m)}(y)\\
    &=\mathit{\Gamma } ( \vartheta )\sum_{n=0}^{\infty}\mathbb{P}\left[ H(t) = n\right]\sum_{m=1}^{\infty}\frac{(\lambda_{\theta}n^{\theta})^m(\zeta )_m}{m!} E_{\mu, \vartheta +m\mu}^{\zeta +m}(-\lambda_{\theta}n^{\theta})S_Y^{(m)}(y)\\
    &= \mathit{\Gamma } ( \vartheta )\sum_{m=1}^{\infty}\frac{\lambda_{\theta}^m(\zeta )_m}{m!}S_Y^{(m)}(y)\sum_{n=0}^{\infty} n^{\theta m}\sum_{l=0}^{\infty}\frac{(\zeta  + m)_{l}(-\lambda_{\theta}n^\theta)^l}{l! \mathit{\Gamma } (\mu (l+m)+ \vartheta )}\mathbb{P}\left[ H(t) = n\right]\\
      &= \mathit{\Gamma } ( \vartheta )\sum_{m=1}^{\infty}\frac{\lambda_{\theta}^m(\zeta )_m}{m!}S_Y^{(m)}(y)\sum_{l=0}^{\infty}\frac{(\zeta  + m)_{l}(-\lambda_{\theta})^l}{l! \mathit{\Gamma } (\mu (l+m)+ \vartheta )}\mathbb{E}\left[(H(t))^{\theta(l+m) }\right].
\end{align*}
With the help of the above two parts, the proof is complete.\end{proof}
\begin{corollary}
    Let the random variables $\{Y_i, i\in \mathbb{N}\}$  are absolutely continuous with  pdf $f_y(\cdot)$, then the probability law of $Z(t),\;t\geq 0$  has the following absolutely continuous and discrete component\\
    (i) the absolutely continuous component is given by
   \begin{equation*}
    h_{Z}(y,t) =  \mathit{\Gamma } ( \vartheta )\sum_{m=1}^{\infty}\frac{\lambda_{\theta}^m(\zeta )_m}{m!}\sum_{l=0}^{\infty}\frac{(\zeta  + m)_{l}(-\lambda_{\theta})^l}{l! \mathit{\Gamma } (\mu (l+m)+ \vartheta )}\mathbb{E}\left[(H(t))^{\theta(l+m) }\right]f_{Y}^{(m)}(y),
\end{equation*}
where $f_{Y}^{(m)}(\cdot)$ is the m-fold convolution of $f_{Y}(\cdot)$.\\
(ii) the discrete component of the probability law is
\begin{align*}
    \mathbb{P}\left[Z(t)=0\right] &= \sum_{n=0}^{\infty}\mathbb{P}\left[ \mathcal{X}_{\mu, \vartheta }^{\zeta ,\theta}(n)= 0\right]\mathbb{P}\left[ H(t) = n\right]\\
    &=  \sum_{n=0}^{\infty}\mathit{\Gamma } ( \vartheta )E_{\mu, \vartheta }^{\zeta }(-\lambda_{\theta}n^{\theta})\mathbb{P}\left[H(t) = n\right],
\end{align*}
where 
\begin{align*}
    \mathbb{P}\left[ \mathcal{X}_{\mu, \vartheta }^{\zeta ,\theta}(t)= 0\right] &=  W_{X}(0,t) -  W_{X}(y,t),\;\;\; y<0\\
    &= \mathit{\Gamma } ( \vartheta )E_{\mu, \vartheta }^{\zeta }(-\lambda_{\theta}t^{\theta}).
\end{align*}
\end{corollary}
  \begin{corollary}\label{cor22}
    If $\{Y_n, n\geq 1\}$ are integer-valued discrete random variables, we set
    \begin{equation*}
        b_s^{*m} = \mathbb{P}\left[Y_1+Y_2+\cdots +Y_m= s\right],\;\; s\in\mathbb{N}_{0}.
    \end{equation*}
    Then, for $s\in \mathbb{N}_{0}$, we have that
    \begin{align}\label{corr23}
        \mathbb{P}\left[Z(t)=s\right]= &\;\mathbb{I}_{\{ s= 0\}} \sum_{k=0}^{\infty}\frac{(\zeta )_k\mathit{\Gamma } ( \vartheta )(-\lambda_{\theta})^k}{k!\mathit{\Gamma } (\mu k+ \vartheta )} \mathbb{E}\left[(H(t))^{\theta k}\right]\\ & + \mathit{\Gamma } ( \vartheta )\sum_{m=1}^{\infty}\frac{\lambda_{\theta}^m(\zeta )_m}{m!}b_s^{*m} \sum_{l=0}^{\infty}\frac{(\zeta  + m)_{l}(-\lambda_{\theta})^l}{l! \mathit{\Gamma } (\mu (l+m)+ \vartheta )}\mathbb{E}\left[(H(t))^{\theta(l+m) }\right].\nonumber
    \end{align}
\end{corollary}  

Next, we present the following special cases for the FGCP at L\'evy times.
\begin{example}
    Let $\{Y_n, n\geq 1\}$ are iid Bernoulli random variables having parameter $ p\in[0,1] $ i.e. $\mathbb{P}\left[Y_n = k\right]= p ^k (1-p)^{1-k}$ for $k \in \{0,1\}$. In accordance to the Corollary \ref{cor22} we have 
    \begin{equation*}
        b_s^{*m} = \binom{m}{s} p^s (1- p)^{m-s}.
    \end{equation*}
    Now, evaluating Eq. (\ref{corr23}) for $s=0$, we obtain
\begin{align*}
        \mathbb{P}\left[Z(t)=0\right]= &\;\sum_{k=0}^{\infty}\frac{(\zeta )_k\mathit{\Gamma } ( \vartheta )(-\lambda_{\theta})^k}{k!\mathit{\Gamma } (\mu k+ \vartheta )} \mathbb{E}\left[(H(t))^{\theta k}\right]\\ & + \mathit{\Gamma } ( \vartheta )\sum_{m=1}^{\infty}\frac{\lambda_{\theta}^m(\zeta )_m}{m!}(1- p)^m \sum_{l=0}^{\infty}\frac{(\zeta  + m)_{l}(-\lambda_{\theta})^l}{l! \mathit{\Gamma } (\mu (l+m)+ \vartheta )}\mathbb{E}\left[(H(t))^{\theta(l+m) }\right].
    \end{align*}
Similarly, for $s \in \mathbb{N}$ Eq. (\ref{corr23}) becomes
\begin{equation*}
     \mathbb{P}\left[Z(t)=s\right]=\mathit{\Gamma } ( \vartheta )\sum_{m=s}^{\infty}\frac{\lambda_{\theta}^m(\zeta )_m}{m!}\binom{m}{s} p^s (1- p)^{m-s} \sum_{l=0}^{\infty}\frac{(\zeta  + m)_{l}(-\lambda_{\theta})^l}{l! \mathit{\Gamma } (\mu (l+m)+ \vartheta )}\mathbb{E}\left[(H(t))^{\theta(l+m) }\right].
\end{equation*}
\end{example}
\begin{example}
    Let $\{Y_n, n\geq 1\}$ are iid Poisson random variables having parameter $\rho>0
    $ i.e.  $\mathbb{P}\left[Y_n = k\right]= \frac{e^{-m\rho}\rho^k}{k!}$ for $k \in \mathbb{N}_0$. In accordance to the Corollary \ref{cor22} we have 
    \begin{equation*}
        b_s^{*m} = \frac{e^{-m\rho}(m\rho)^s}{s!},\;\;s \in \mathbb{N}_0.
    \end{equation*}
Now, evaluating Eq. (\ref{corr23}) for  $s \in \mathbb{N}_0$ we get
\begin{align*}
        \mathbb{P}\left[Z(t)=s\right]= &\;\mathbb{I}_{\{ s= 0\}} \sum_{k=0}^{\infty}\frac{(\zeta )_k\mathit{\Gamma } ( \vartheta )(-\lambda_{\theta})^k}{k!\mathit{\Gamma } (\mu k+ \vartheta )} \mathbb{E}\left[(H(t))^{\theta k}\right]\\ & + \mathit{\Gamma } ( \vartheta )\sum_{m=1}^{\infty}\frac{\lambda_{\theta}^m(\zeta )_m}{m!}\frac{e^{-m\rho}(m\rho)^s}{s!} \sum_{l=0}^{\infty}\frac{(\zeta  + m)_{l}(-\lambda_{\theta})^l}{l! \mathit{\Gamma } (\mu (l+m)+ \vartheta )}\mathbb{E}\left[(H(t))^{\theta(l+m) }\right].
    \end{align*}
At $s=0$, it reduces to
\begin{align*}
        \mathbb{P}\left[Z(t)=0\right]= &\;\sum_{k=0}^{\infty}\frac{(\zeta )_k\mathit{\Gamma } ( \vartheta )(-\lambda_{\theta})^k}{k!\mathit{\Gamma } (\mu k+ \vartheta )} \mathbb{E}\left[(H(t))^{\theta k}\right]\\ & + \mathit{\Gamma } ( \vartheta )\sum_{m=1}^{\infty}\frac{e^{-m\rho}\lambda_{\theta}^m(\zeta )_m}{m!} \sum_{l=0}^{\infty}\frac{(\zeta  + m)_{l}(-\lambda_{\theta})^l}{l! \mathit{\Gamma } (\mu (l+m)+ \vartheta )}\mathbb{E}\left[(H(t))^{\theta(l+m) }\right].
    \end{align*}
\end{example}

\begin{example}
      Let $\{Y_n, n\geq 1\}$ are iid geometrically random variables having parameter $ p\in[0,1] $ i.e. $\mathbb{P}\left[Y_n = k\right]= p (1-p)^{k-1}$ for $k \in \mathbb{N}$. In accordance to the Corollary \ref{cor22}, we have 
    \begin{equation*}
        b_s^{*m} = \binom{s-1}{m-1} p^m (1- p)^{s-m}, \;\;\;\;s\geq m.
    \end{equation*}
    As we know, the sum of independent random variables having geometric distribution with the same parameter follows the negative binomial distribution.  For $s \in \mathbb{N}_0$ evaluate Eq. (\ref{corr23})
 \begin{align*}
        \mathbb{P}\left[Z(t)=s\right]= &\;\mathbb{I}_{\{ s= 0\}} \sum_{k=0}^{\infty}\frac{(\zeta )_k\mathit{\Gamma } ( \vartheta )(-\lambda_{\theta})^k}{k!\mathit{\Gamma } (\mu k+ \vartheta )} \mathbb{E}\left[(H(t))^{\theta k}\right]\\ & + \mathit{\Gamma } ( \vartheta )\sum_{m=1}^{s}\frac{\lambda_{\theta}^m(\zeta )_m}{m!}\binom{s-1}{m-1} p^m (1- p)^{s-m} \sum_{l=0}^{\infty}\frac{(\zeta  + m)_{l}(-\lambda_{\theta})^l}{l! \mathit{\Gamma } (\mu (l+m)+ \vartheta )}\mathbb{E}\left[(H(t))^{\theta(l+m) }\right].
    \end{align*}
\end{example}

\section{Some interconnections and applications}\label{sec:appl}
 It is well known that the Bell polynomial plays a crucial role in the various disciplines of applied sciences. The Bell polynomials were reported in terms of the recurrence relation of the Bell numbers by Gould and Quaintance \cite{boubellouta2020some}. It is also known that the Bell polynomials are closely associated with the  $n$-th order moments of the
 Poisson random variables. Later on, the fractional generalization of the Bell polynomials was introduced and studied in \cite{laskin2009some}. Recently, Laskin \cite{laskin2024new} introduced the generalized fractional Bell polynomials (GFBP) denoted by $\mathcal{B}_{G}(x,m)$ by utilizing the pmf of the FCP.  The GFBP of $m$-th order is given by 
    \begin{equation*}
        \mathcal{B}_{G}(x,m) = \mathit{\Gamma } ( \vartheta )\sum_{n=0}^{\infty}n^m (\zeta )_n\frac{ x^n }{n!}\sum_{k=0}^{\infty}\frac{(\zeta +n)_{k}}{k!\mathit{\Gamma } (\mu(n+k)+ \vartheta )}(-x)^k, \;\;\; m\in \mathbb{N}_0.
    \end{equation*}
 One can also find some recent generalizations of the Bell polynomials in \cite{soni2024probabilistic, kim2023probabilistic}. In this section, we first define the subordinated version of the GFBP and subsequently discover an interconnection between the moments of the TCFCP with the subordinated generalized fractional Bell polynomial.
\subsection{Subordinated generalized  fractional Bell polynomials and subordinated generalized fractional Bell numbers}
With the help of Eq. (\ref{pmfnew}), we introduce a L\'{e}vy subordinated version of the GFBP called subordinated generalized fractional Bell polynomials (SGFBP). For $m \in \mathbb{N}_0$ we denote the $m$-th order SGFBP by 
\begin{align*}
    \mathcal{B}_{SG}(x,m)&= \sum_{n=0}^{\infty}n^mz(n,x),
    \end{align*}
   where the expression for $z(n,x)$ is given in Eq. (\ref{pd}). Hence, we get the series representation as
\begin{align}\label{bell}
    \mathcal{B}_{SG}(x,m)&= \mathit{\Gamma } ( \vartheta )\sum_{n=0}^{\infty}\frac{n^m (\zeta )_n }{n!}\sum_{k=0}^{\infty}\frac{(-1)^k(\zeta +n)_{k}}{k!\mathit{\Gamma } (\mu(n+k)+ \vartheta )}\mathbb{E}\left[X^{\theta(n+k)}\right] ,\;\;\;\mathcal{B}_{SG}(x,0)=1.
\end{align}
\noindent When we take $x=1$ in the above equation, we get the Bell numbers, which we introduce as the  subordinated generalized fractional Bell numbers (SGFBN) and are denoted by $\mathcal{B}_{SG}(m)$
\begin{align*}
    \mathcal{B}_{SG}(m)= \mathcal{B}_{SG}(x,m)\big{|}_{x=1} &= \mathit{\Gamma } ( \vartheta )\sum_{n=0}^{\infty}\frac{n^m (\zeta )_n }{n!}\sum_{k=0}^{\infty}\frac{(-1)^k(\zeta +n)_{k}}{k!\mathit{\Gamma } (\mu(n+k)+ \vartheta )}\\
    &=\mathit{\Gamma } ( \vartheta )\sum_{n=0}^{\infty}\frac{n^m }{n!}E_{\mu, \vartheta }^{\zeta (n)}(-1).\nonumber
\end{align*}
\begin{remark}
    Relationship between the SGFBP and the GFBP.\\
\noindent(i) Let us consider the random variable $X$ to be an $\alpha$-stable subordinator $S_{\alpha}(t),\;t\geq0$ and $\alpha\in(0,1) $ having the fractional moments of order $p>0$ is given as follows (see \cite{beghin2020tempered})
    \begin{equation*}
        \mathbb{E}\left[S_{\alpha}(t)^p\right] = \frac{\mathit{\Gamma } ( 1-p/\alpha )}{\mathit{\Gamma } ( 1-p )}t^{p/\alpha}.
    \end{equation*}
    Substituting it in Eq. (\ref{bell}) we obtain the following expression
    \begin{align*}
        \mathcal{B}_{SG}(t,m)&= \mathit{\Gamma } ( \vartheta )\sum_{n=0}^{\infty}\frac{n^m (\zeta )_n }{n!}\sum_{k=0}^{\infty}\frac{(-1)^k(\zeta +n)_{k}}{k!\mathit{\Gamma } (\mu(n+k)+ \vartheta )}\frac{\mathit{\Gamma } ( 1-\frac{\theta(n+k)}{\alpha} )}{\mathit{\Gamma } ( 1-{\theta(n+k)} )}t^{\frac{\theta(n+k)}{\alpha}}.
    \end{align*}
    Considering the limiting case for  $\alpha=1 $ and subsequently replacing $t^{\theta}$ by $x$, we get
    \begin{equation*}
        \mathcal{B}_{SG}(x,m) = \mathcal{B}_{G}(x,m).
    \end{equation*}
    \noindent(ii) Now, let us suppose that random variable $X$  be an incomplete gamma subordinator ${D_{\alpha}(t)},\;t\geq0$ for $\alpha\in(0,1]$ having the fractional moments of order $p$ for $p\leq \alpha$ has asymptotic behavior as (see \cite{beghin2021levy})
    \begin{equation*}
        \mathbb{E}\left[D_{\alpha}(t)^p\right] \simeq \frac{\mathit{\Gamma } ( 1-p/\alpha )}{\mathit{\Gamma } ( 1-p )}t^{p/\alpha},\;\;\;\; t \to \infty.
    \end{equation*}
   Substituting this in Eq. (\ref{bell}), we obtain the following expression
    \begin{align*}
       \mathcal{B}_{SG}(t,m)& \simeq \mathit{\Gamma } ( \vartheta )\sum_{n=0}^{\infty}\frac{n^m (\zeta )_n }{n!}\sum_{k=0}^{\infty}\frac{(-1)^k(\zeta +n)_{k}}{k!\mathit{\Gamma } (\mu(n+k)+ \vartheta )}\frac{\mathit{\Gamma } ( 1-\frac{\theta(n+k)}{\alpha} )}{\mathit{\Gamma } ( 1-{\theta(n+k)} )}t^{\frac{\theta(n+k)}{\alpha}}.
    \end{align*}
    Let $\alpha=1 $ and by replacing $t^{\theta}$ by $x$, we get
    \begin{equation*}
        \mathcal{B}_{SG}(x,m) \simeq \mathcal{B}_{G}(x,m),
    \end{equation*}
  
   \noindent(iii) Let us consider the random variable $X$ to be a tempered stable subordinator $S_{\alpha, \varphi }(t),\;t\geq0$ and $\alpha\in(0,1),\;  \varphi >0 $ having the fractional moments of order $p$  has asymptotic behaviour as (see \cite{kumar2019fractional}) 
    \begin{equation*}
        \mathbb{E}\left[S_{\alpha, \varphi }(t)^p\right] \sim \left(\alpha \varphi^{\alpha-1} t\right)^{p}, \;\;\text{as} \;t\rightarrow\infty.
    \end{equation*}
    Substituting this in Eq. (\ref{bell}) we obtain the following expression
    \begin{align*}
       \mathcal{B}_{SG}(t,m)& \sim \mathit{\Gamma } ( \vartheta )\sum_{n=0}^{\infty}\frac{n^m (\zeta )_n }{n!}\sum_{k=0}^{\infty}\frac{(-1)^k(\zeta +n)_{k}}{k!\mathit{\Gamma } (\mu(n+k)+ \vartheta )}\left(\alpha \varphi^{\alpha-1} t\right)^{\theta(n+k)}.
    \end{align*}
   Considering the limiting case for $\alpha=1 $ and subsequently replacing $t^{\theta}$ by $x$, we get
    \begin{equation*}
       \mathcal{B}_{SG}(x,m) \sim \mathcal{B}_{G}(x,m).
    \end{equation*}
    \end{remark}
  
\begin{remark}
    It is easy to observe that the SGFBN and the generalized fractional Bell number (GFBN) are equal 
    (see \cite{laskin2024new})
    \begin{equation*}
          \mathcal{B}_{SG}(m)=\mathcal{B}_{G}(1,m).
    \end{equation*}
\end{remark}
\subsection{Generating function of the SGFBP}
The generating function of the SGFBP, denoted by $\mathcal{F}_{SG}(s,x)$, can be evaluated as 
\begin{align}\label{gf}
    \mathcal{F}_{SG}(s,x)=\sum_{m=0}^{\infty}\frac{s^m}{m!} \mathcal{B}_{SG}(x,m).
\end{align}
By differentiating $\mathcal{F}_{SG}(s,x)$ 
$m$ times with respect to $s$ and evaluating it at $s=0$, we obtain $\mathcal{B}_{SG}(x,m)$ as 
\begin{equation*}
    \mathcal{B}_{SG}(x,m)= \frac{\partial^m}{\partial s^m}\mathcal{F}_{SG}(s,x)\big{|}_{s=0}.
\end{equation*}
Substituting Eq. (\ref{bell}) in Eq. (\ref{gf}), we obtain $\mathcal{F}_{SG}(s,x)$ 
\begin{align*}
     \mathcal{F}_{SG}(s,x)=\sum_{m=0}^{\infty}\frac{s^m}{m!}  \mathit{\Gamma } ( \vartheta )\sum_{n=0}^{\infty}\frac{n^m (\zeta )_n }{n!}\sum_{k=0}^{\infty}\frac{(-1)^k(\zeta +n)_{k}}{k!\mathit{\Gamma } (\mu(n+k)+ \vartheta )}\mathbb{E}\left[X^{\theta(n+k)}\right].
\end{align*}
Now, by evaluating the sum over $m$, we get
\begin{align*}
     \mathcal{F}_{SG}(s,x)=& \; \mathit{\Gamma } ( \vartheta )\sum_{n=0}^{\infty}\frac{e^{sn} (\zeta )_n }{n!}\sum_{k=0}^{\infty}\frac{(-1)^k(\zeta +n)_{k}}{k!\mathit{\Gamma } (\mu(n+k)+ \vartheta )}\mathbb{E}\left[X^{\theta(n+k)}\right]\\
     =&\;\sum_{n=0}^{\infty}e^{sn} z(n,x).
\end{align*}
With the help of Eq. (\ref{mgsx}), we get
\begin{align}\label{fsx}
    \mathcal{F}_{SG}(s,x)=  \mathit{\Gamma } ( \vartheta )\sum_{l=0}^{\infty}\frac{(\zeta )_l}{l!\mathit{\Gamma } (\mu l+ \vartheta )}\mathbb{E}\left[(X^{\theta}(e^{s}-1))^l\right].
\end{align}
Put $x=1$ in Eq. (\ref{fsx}) to obtain the generating function of the SGFBN $\mathcal{B}_{SG}(m)$ 
\begin{equation*}
    \mathcal{F}_{SG}(s)=\mathcal{F}_{SG}(s,x)\big{|}_{x=1}=\mathit{\Gamma } ( \vartheta )\sum_{l=0}^{\infty}\frac{(\zeta )_l}{l!\mathit{\Gamma } (\mu l+ \vartheta )}(e^{s}-1)^l.
\end{equation*}
Therefore, we get 
\begin{equation*}
    \mathcal{F}_{SG}(s)= \mathcal{F}_{G}(s),
\end{equation*}
where $\mathcal{F}_{G}(s)$ denotes the generating function of the GFBN (see \cite{laskin2024new}).

\subsection{Moments of TCFCP}
The SGFBP is now used to evaluate the moments of the process TCFCP $Z(t,\lambda _{\theta }),\;t\geq 0$. By the definition of the $p^{\text{th}}$ order moment, we have that
\begin{equation*}
    \mathbb{E}\{[Z(t,\lambda _{\theta })]^p\} = \sum_{n=0}^{\infty} n^p z(n,t). 
\end{equation*}
It can be easily followed that $\mathbb{E}\{[Z(t,\lambda _{\theta })]^p\}$ can be expressed in terms of the SGFBP as
\begin{equation*}
     \mathbb{E}\{[Z(t,\lambda _{\theta })]^p\} =   \mathcal{B}_{SG}(\lambda_{\theta}(H(t))^{\theta},m).
\end{equation*}
In the next subsection, we present an application of the FCP in the shock deterioration model.
\subsection{Application of the FCP in shock deterioration model}
%Here, we present a FGCP (see \cite{laskin2024new}) based shock model. 
 Civil structures undergo several natural and man-made risks, such as heavy automobiles, strong winds, earthquakes, and
 flooding. These damages physically harm the in-service civil structures, thereby reducing their safety below a reasonable threshold thus affecting their performance (see \cite{wang2022explicit,obrien2014lifetime,liu2020temperature,lee2016new,alos2014analysis}).  The resistance degrades due to repeated attacks. Accurate prediction of the serviceability of affected structures is difficult due to uncertainties in both their structural properties and the external stresses they encounter.  Within this framework, the structural reliability evaluation serves as a valuable instrument for quantifying the safety level of deteriorated structures subjected to various threats, producing probability-based recommendations for structural design and maintenance decisions. Recently, Wang \cite{wang2022explicit} introduced a compound Poisson process-based shock deterioration model for the reliability assessment of aging structures. We extend the model proposed in \cite{wang2022explicit} by utilizing the FGCP in the model. The FCP incorporates memory effects, which are essential for modelling systems where past shocks influence future events.  It also allows flexible interarrival times, unlike the Poisson process, making it better suited for complex real-world phenomena.
The FGCP is a stochastic process in which the FCP counts the happening of $Y_i$ in time $t$ is defined as 
\begin{equation}\label{shock}
	\mathcal{X}_{\mu, \vartheta }^{\zeta ,\theta}(t) =  \sum_{i=0}^{\mathcal{N}_{\mu, \vartheta }^{\zeta ,\theta}(t)}Y_i, t\geq 0
\end{equation}
with $Y_0 \equiv 0$. Taking into account the degradation of aging structure due to shock, we can use $\mathcal{X}_{\mu, \vartheta }^{\zeta ,\theta}(t)$ to determine the difference in the initial resistance and the deteriorated resistance at time $t$.\\
Next, the pdf of $\mathcal{X}_{\mu, \vartheta }^{\zeta ,\theta}(t)$ denoted by $f_{\mathcal{X}}(y,t)$, is evaluated with the help of law of total probability.
\begin{align}\label{nfold}
	f_{\mathcal{X}}(y,t) =& \sum_{n=1}^{\infty}\mathcal{P}_{\mu, \vartheta }^{\zeta ,\theta}(n,t) f_Y(y)\ast f_Y(y)\ast \cdots \ast f_Y(y)\nonumber\\
	=& \sum_{n=1}^{\infty}\mathcal{P}_{\mu, \vartheta }^{\zeta ,\theta}(n,t) f_Y^{<n>}(y),
\end{align}
where $f_Y(y)$ denotes the pdf of $Y_k$ and $f_Y^{<n>}(y)$ represents the $n$-fold convolution of $f_Y(y)$. The widely used model for the gradual degradation of resistance of the structure is the gamma process, which means that a gamma-distributed random variable is used to simulate the amount of each shock deterioration increment (see \cite{wang2022explicit,wang2015realistic,lawless2004covariates}). Now, if $Y_k$'s are iid gamma distributed, the pdfs for $Y_k$, and its convolution is given by
\begin{align*}
	f_Y(y) &= \frac{\left(y/b\right)^{a-1}e^{-y/b}}{b\mathit{\Gamma } (a)},\; y\geq0,\\
	f_Y^{<n>}(y) &= \frac{\left(y/b\right)^{an-1}e^{-y/b}}{b\mathit{\Gamma } (an)},\; y\geq0.
\end{align*}
Hence, the Eq. (\ref{nfold}) becomes
\begin{align}\label{pdfsh}
	f_{\mathcal{X}}(y,t) =& \sum_{n=1}^{\infty}\frac{(\zeta )_n \mathit{\Gamma } ( \vartheta )(\lambda_{\theta}t^{\theta})^n}{n!}  \frac{\left(y/b\right)^{an-1}e^{-y/b}}{b\mathit{\Gamma } (an)}E_{\mu, \vartheta+n\mu }^{\zeta+n }(-\lambda_{\sigma}t^{\sigma})\nonumber\\
	=&\; \frac{\mathit{\Gamma } ( \vartheta ) e^{-y/b}}{y}\sum_{n=1}^{\infty}\frac{(\zeta )_n (\lambda_{\theta}t^{\theta}(y/b)^a)^n}{n!\mathit{\Gamma } (an)}E_{\mu, \vartheta+n\mu }^{\zeta+n }(-\lambda_{\sigma}t^{\sigma}).
\end{align}
We define the $ \Upsilon (p,q,s)$ function as follows
\begin{equation*}
	\Upsilon (p,q,s) = \sum_{n=1}^{\infty} \frac{(\zeta )_n (pq^s)^n}{n!\mathit{\Gamma } (sn)}E_{\mu, \vartheta+n\mu }^{\zeta+n }(-p).
\end{equation*}
With this the Eq. (\ref{pdfsh}) becomes
\begin{equation}\label{pll}
	f_{\mathcal{X}}(y,t) = \frac{\mathit{\Gamma } ( \vartheta ) e^{-y/b}}{y} \Upsilon (\lambda_{\sigma}t^{\sigma},y/b,a).
\end{equation}
Hence, the cdf of $\mathcal{X}_{\mu, \vartheta }^{\zeta ,\theta}(t)$ can be obtained by integrating Eq. (\ref{pll})
\begin{align*}
	F_{\mathcal{X}}(y,t)=&\; \mathit{\Gamma } ( \vartheta )E_{\mu, \vartheta }^{\zeta }(-\lambda_{\theta}t^{\theta})+\int_{0}^{y}f_{\mathcal{X}}(x,t)\mathrm{d} x\\ 
	=&\; \mathit{\Gamma } ( \vartheta )E_{\mu, \vartheta }^{\zeta }(-\lambda_{\theta}t^{\theta})+\int_{0}^{y}\frac{\mathit{\Gamma } ( \vartheta ) e^{-x/b}}{x} \Upsilon (\lambda_{\sigma}t^{\sigma},x/b,a)\mathrm{d} x.
\end{align*}
Taking into account the probability of $\mathcal{N}_{\mu, \vartheta }^{\zeta ,\theta}(t)=0$ the term $ \mathit{\Gamma } ( \vartheta )E_{\mu, \vartheta }^{\zeta }(-\lambda_{\theta}t^{\theta})$ has been included.
Next, we show that the $ \Upsilon (p,q,s)$ is bounded above. As  $\mathit{\Gamma } (x) > 0.8856$ for $x>0$ (see \cite{deming1935minimum}) and with the help of the Taylor series, we get 
\begin{align*}
	\Upsilon (p,q,s) &< \frac{1}{0.8856}\sum_{n=1}^{\infty} \frac{(\zeta )_n (pq^s)^n}{n!}E_{\mu, \vartheta }^{\zeta(n) }(-p)\\
	&= \frac{1}{0.8856}\left[E_{\mu, \vartheta }^{\zeta}(p(q^s-1))-1\right],
\end{align*}
which gives an upper bound of $ \Upsilon (p,q,s)$.\\
Next, we discuss the reliability assessment of the proposed model.
Several criteria can be utilized to determine whether a structure fails or remains operational. Here, we consider the criterion that the instance at which the cumulative resistance deterioration exceeds the pre-decided threshold indicates the damage dominant failure
 mode. A coefficient based on the desired reliability level can be multiplied by the total load effects to define the threshold. The linear combination of the gradual and shock deterioration represents the total resistance degradation (see \cite{iervolino2013gamma, wang2017reliability}). The primary reason
 for gradual deterioration is environmental factors, which are influenced by structural
 materials and, in turn, by the dominant deterioration mechanisms.  Conversely, severe attacks are usually the source of shock deterioration. The resistance $\mathcal{R}(t)$ at time $t$ is calculated mathematically as
\begin{equation*}
	\mathcal{R}(t) = \mathcal{R}_0 - \mathcal{X}_{\mu, \vartheta }^{\zeta ,\theta}(t) - \mathcal{S}(t),
\end{equation*}
where $\mathcal{R}_0$ represents the initial resistance, $\mathcal{X}_{\mu, \vartheta }^{\zeta ,\theta}(t)$ and $\mathcal{S}(t)$ 
represents the gradual and shock resistance degradation up to time $t$, respectively. Let us consider a time period $[0,\tau]$ in which the structure is accounted as survived if the deteriorated resistance at time $\tau$, i.e. $\mathcal{R}(\tau)$, is not lesser than the predefined threshold $k_p$, that is
\begin{equation*}
	Y(0,\tau)= \mathbb{P}\left[\mathcal{R}(\tau)\geq k_p \right]= \mathbb{P}\left[\mathcal{R}_0 - \mathcal{X}_{\mu, \vartheta }^{\zeta ,\theta}(\tau) - \mathcal{S}(\tau)\geq k_p \right],
\end{equation*}
where $ Y(0,\tau)$ denotes the survival probability in time $[0,\tau]$. With the help of the law of total probability, we have 
\begin{align}\label{sur}
	Y(0,\tau) =&\; \mathbb{P}\left[ \mathcal{X}_{\mu, \vartheta }^{\zeta ,\theta}(\tau) \leq \mathcal{R}_0-\mathcal{S}(\tau) - k_p \right]\nonumber\\
	=& \int _{\mathcal{R}_0}\int _{\mathcal{S}(\tau)}F_{\mathcal{X}}(r_0-s-k_p,\tau)f_{\mathcal{S}(\tau)}(s)f_{\mathcal{R}_0}(r_0)\mathrm{d} s\mathrm{d} r_0,
\end{align}
where $f_{\mathcal{S}(\tau)}(s)$ and $f_{\mathcal{R}_0}(r_0)$ are the pdfs of  $\mathcal{S}(\tau)$ and $\mathcal{R}_0$ respectively. In particular, if the initial resistance has negligible uncertainty (given by $r_0$), as is generally the case in practical engineering (see \cite{ellingwood1985probabilistic}),  Eq. (\ref{sur}) becomes
\begin{equation*}
	Y(0,\tau) = \int _{\mathcal{S}(\tau)}F_{\mathcal{X}}(r_0-s-k_p,\tau)f_{\mathcal{S}(\tau)}(s)\mathrm{d} s.
\end{equation*}
The probability of failure of structure over time interval $[0,\tau]$, represented by $Q(0,\tau)$ can be calculated as 
\begin{equation*}
	Q(0,\tau) = 1- Y(0,\tau).
\end{equation*}
Let $T$ be the time of failure, and its cdf be $F_T$. Therefore, by definition, we get
\begin{equation*}
	F_T(T) = \mathbb{P}\left[T\leq \tau\right] = Q(0,\tau) = 1- Y(0,\tau).
\end{equation*}
It indicates that once the time-dependent reliability is established, the cdf of the time to failure can be calculated.
\begin{remark}
	When $\mu= \vartheta= \zeta =\theta = 1$, the FGCP based shock deterioration model Eq. (\ref{shock}) reduces to the compound Poisson based shock deterioration model as done in \cite{wang2022explicit}.
\end{remark}

%\bibliography{referencee22nov}
\bibliographystyle{abbrv}

%\vspace{-1cm}
\end{document}